\documentclass[review,12pt]{elsarticle}
\usepackage{amsmath}
\usepackage{amsfonts}
\usepackage{geometry}   
\usepackage[parfill]{parskip}    
\usepackage{graphicx}
\usepackage{amssymb}
\usepackage{epstopdf}
\usepackage{psfrag}
\newcommand{\figref}[1]{{Figure~\ref{#1}}}
\newcommand {\bea}{\begin{eqnarray}}
\newcommand {\ea}{\end{eqnarray}}
\newtheorem{proposition}{Proposition}[section]
\newtheorem{theorem}{Theorem}[section]
\newtheorem{Assumption}{Assumption}[section]
\newtheorem{lemma}{Lemma}[section]
\newtheorem{remark}{Remark}[section]

\newtheorem{corollary}{Corollary}[section]

\newenvironment{proof}[1][Proof]{\textbf{#1.} }{\hspace{\stretch{1}}\rule{0.5em}{0.5em}}

\usepackage{xspace}

\newcommand{\Dt}{\Delta t}
\usepackage{subfigure}
\usepackage{color}

\newcommand{\thmref}[1]{{Theorem~\ref{#1}}}
\newcommand{\lemref}[1]{{Lemma~\ref{#1}}}
\newcommand{\secref}[1]{{Section~\ref{#1}}}
\newcommand{\assref}[1]{{Assumption~\ref{#1}}}
\newcommand{\propref}[1]{{Proposition~\ref{#1}}}
\newcommand{\coref}[1]{{Corollary~\ref{#1}}}

\journal{Computers \& Mathematics with Applications}
\begin{document}
\begin{frontmatter}
\title{Optimal strong convergence rates of numerical methods for  semilinear parabolic SPDE driven by Gaussian noise and Poisson random measure}

\author[jdm]{Jean Daniel Mukam}
\ead{jean.d.mukam@aims-senegal.org}
\address[jdm]{Fakult\"{a}t f\"{u}r Mathematik, Technische Universit\"{a}t Chemnitz, 09126 Chemnitz, Germany}

\author[ata,atb,atc]{Antoine Tambue}
\cortext[cor1]{Corresponding author}
\ead{antonio@aims.ac.za}
\address[ata]{Department of Computing Mathematics and Physics,  Western Norway University of Applied Sciences, Inndalsveien 28, 5063 Bergen}
\address[atc]{Center for Research in Computational and Applied Mechanics (CERECAM), and Department of Mathematics and Applied Mathematics, University of Cape Town, 7701 Rondebosch, South Africa.}
\address[atb]{The African Institute for Mathematical Sciences(AIMS) of South Africa,
6-8 Melrose Road, Muizenberg 7945, South Africa}

%

\begin{abstract}
This paper deals with the numerical approximation of semilinear parabolic stochastic partial differential equation (SPDE) driven simultaneously by Gaussian noise and Poisson random measure, more realistic in modeling real world phenomena.
  The SPDE is discretized in space with the standard finite element method and in time with the linear implicit Euler method or  an exponential integrator, more efficient and stable for stiff problems.
   We prove the strong convergence of the fully discrete schemes toward the mild solution. The results reveal how convergence orders depend on the regularity of the noise and the initial data. 
   In addition, we exceed the classical orders $1/2$ in time and $1$ in space achieved in the literature when dealing with SPDE driven by Poisson measure with less regularity assumptions on  the nonlinear drift function. 
   In particular, for trace class multiplicative Gaussian noise we achieve convergence order $\mathcal{O}(h^2+\Delta t^{1/2})$. 
   For additive trace class Gaussian noise and an appropriate jump function, we achieve convergence order $\mathcal{O}(h^2+\Delta t)$. 
    Numerical experiments to sustain the theoretical results are provided.
\end{abstract}

\end{frontmatter}
\section{Introduction}
\label{intro}
We consider the SPDE defined in $\Lambda\subset \mathbb{R}^d$, $d=1,2,3$,  with initial value  of the following type
\begin{eqnarray}
\label{model}
dX(t)=[AX(t)+F(X(t))]dt+B(X(t))dW(t)+\int_{\chi}G(z, X(t))\widetilde{N}(dz,dt), \quad X(0)=X_0, 
\end{eqnarray}
 $  t\in(0,T]$ on the Hilbert space $H=\left(L^2(\Lambda), \langle.,.\rangle_H, \Vert.\Vert\right)$, where  $T>0$ is the final time,   
$A$ is a linear operator which is unbounded,  not necessarily self-adjoint and is assumed to generate an analytic  semigroup $S(t)$.
The noise  $W(t)=W(x,t)$ is a $H$-valued $Q$-Wiener process defined in a filtered probability space $(\Omega,\mathcal{F}, \mathbb{P}, \{\mathcal{F}_t\}_{t\geq 0})$, where $Q : H\longrightarrow H$ is a positive and linear self-adjoint operator. The filtration is assumed to fulfill the usual assumptions (see e.g., \cite[Definition 2.1.11]{Prevot}). The $Q$-Wiener process $W(t)$ can be represented as follows \cite{Prevot, Prato}
\begin{eqnarray}
W(x, t)=\sum_{i=0}^{\infty}\sqrt{q_i}e_i(x)\beta_i(t),
\end{eqnarray}
where $(q_i)_{i=0}^{\infty}$ are the eigenvalues of $Q$ and $(e_i)_{i=0}^{\infty}$ the corresponding eingenfunctions, and $\beta_i$ are the independent and identically distributed standard normal distributions. 
The mark set $\chi$ is defined by $\chi:=H\setminus\{0\}$.  For a given  set  $\Gamma$,  we denote by $\mathcal{B}(\Gamma)$ the smallest $\sigma$-algebra containing all open sets of $\Gamma$. Let  $(\chi, \mathcal{B}(\chi), \nu)$ be a $\sigma$-finite measurable space and $\nu$ ( with $\nu\not\equiv 0$) a L\'{e}vy measure on $\mathcal{B}(\chi)$ such that
\begin{eqnarray}
\nu(\{0\})=0\quad \text{and}\quad \int_{\chi}\min(\Vert z\Vert^2, 1)\nu(dz)<\infty.
\end{eqnarray}
 Let $N(dz, dt)$ be the $H$-valued  Poisson  distributed $\sigma$-finite measure on the product $\sigma$-algebra $\mathcal{B}(\chi)$ and $\mathcal{B}(\mathbb{R}_{+})$ with intensity $\nu(dz)dt$, where $dt$ is the Lebesgue measure on $\mathcal{B}(\mathbb{R}_{+})$. In our model problem \eqref{model}, $\widetilde{N}(dz, dt)$ stands for the compensated Poisson random measure defined by
 \begin{eqnarray}
 \widetilde{N}(dz, dt):=N(dz, dt)-\nu(dz)dt.
 \end{eqnarray}
 Note that $\widetilde{N}(dz, dt)$ is a noncontinuous martingale with mean $0$, see e.g., \cite{Rudiger}. The Wiener process $W$ and the compensated Poisson measure $\widetilde{N}$ are supposed to be independent. Precise assumptions on the nonlinear functions $F$, $B$ and $G$ to ensure the existence of the mild solution of \eqref{model} will be given in the following section.
Equations of type \eqref{model} with only Gaussian noise   are used to model many real world phenomena in different fields such as   biology, chemistry, physics, see e.g.,  \cite{ATthesis,Shardlow,SebaGatam}. But for many phenomenon, modeling with only Gaussian noise is unsatisfactory and  less realistic. For instance, in finance, the unpredictable nature of many events such as market crashes, announcements made by the central banks,  changing    credit risk, insurance in a changing risk, changing face of operational risk \cite{Tankov, Platen1} might have sudden and significant impacts on the stock price. For such phenomena, to obtain a more realistic model equation, it is recommended to incorporate non-Gaussian noise such as L\'{e}vy process or Poisson random measure. For  more applications   with non-Gaussian noise, see  for example \cite[Chapter 3]{Walsh} for the model of the cable equation in neurophysiology, \cite{Levy} for phenomena in fluid mechanics, polymer chemistry and economic science. 
        Since explicit solutions of \eqref{model}  do not  exist, numerical algorithms are  good tools to provide realistic approximations. Numerical methods for SPDE driven only by Gaussian noise are widely investigated in the scientific literature (see e.g., \cite{Jentzen1,Jentzen2,Gaby1,Kovac1,Raphael,Antonio1, AntonioRev1,AntonioRev2,Xiaojie2,Xiaojie3,Yan2} and references therein). The case of time-fractional SPDEs driven by Gaussian noise have recently received some attentions (see e.g. \cite{Zou1,Zou2} and references therein). The case of SPDE driven by fractional Brownian motion was also recently investigated, see e.g. \cite{XiaojieRev1}.  If we turn our attention to the case of SPDE driven only by non-Gaussian noise, the list of references become remarkably short (see e.g., \cite{Barth1,Barth2,Barth3,Erika3,Erika2,Erika1}). 
        In finite dimension, numerical schemes for  ordinary stochastic  differential equations (SDE) driven by Brownian motion and Poisson process have been extensively studied (see e.g., \cite{Bruti1,Darelotis1,Kloeden1,Kumar,Platen1,Gan} and references therein).  In  infinite dimension, numerical methods for SPDE  \eqref{model} driven by both Gaussian and non-Gaussian noises is quite limited.  
         Recently the finite element method combined with the linear implicit method was applied to \eqref{model} in \cite{Yang} with  a  self adjoint operator $A$.
         The  accuracy of semi-implicit method is limited as the semigroup is approximated by  a rational function of the linear operator (one  of  its  resolvent).
          It is therefore  important to investigate numerical  algorithms  based on exponential integrators where  the discrete semigroup can be evaluated almost exact  for  a given tolerance precision.
          The  huge interest in such numerical integrations is due to the availability since late 1950s of several numerical techniques and  software  for computing the  exponential of matrices such as Krylov subspace technique and Fast L\'{e}ja points  see e.g.,  \cite{Kry,Leja1,Sidje,Leja2}.  In this paper,  we  extend   a stochastic  exponential integrator  provided in \cite{Antonio1,AntonioRev2}  for only Gaussian  process  to   SPDE \eqref{model}.
           Let us recall that the strong convergence orders achieved in \cite{Yang}  for semi  implicit scheme was sub-optimal  of type $\mathcal{O}\left(h^{1-\epsilon}+\Delta t^{1/2-\epsilon}\right)$, where $\epsilon>0$ is small enough for both  multiplicative and additive Wiener and Poisson processes. 
           In fact,  in finite dimension,  for ordinary SDE  driven by multiplicative Wiener  and Poisson processes, the optimal  order of the strong convergence is $1/2$ (see e.g., \cite{Gan,Kloeden1,Platen1,Bruti1}), 
           whereas when dealing with ordinary SDE driven by additive Wiener and Poisson processes the optimal strong convergence order is $1$ (see e.g., \cite{Platen1,Bruti1}).
           In this paper, we extend the result to infinite dimension case  and achieve  the  optimal  strong convergence orders $\mathcal{O}\left(h^2+\Delta t^{1/2}\right)$ for   multiplicative trace class Gaussian noise and  Poisson measure, and  $\mathcal{O}(h^2+\Delta t)$ for additive Gaussian noise and Poisson measure for  our exponential integrator  scheme on SPDE \eqref{model}.  Furthermore,  we   investigate   the semi-implicit Euler method for SPDE \eqref{model} with additive Gaussian and Poisson measure with not necessary self adjoint operator, and  provide appropriate conditions under which the scheme achieves strong convergence orders $\mathcal{O}\left(h^2+\Delta t^{1-\epsilon}\right)$. 
 Note that such optimal convergence orders were also achieved in \cite{Xiaojie2,Xiaojie3,Antjd2,AntonioRev1,AntonioRev2} for SPDE driven only by Gaussian noise.  The case of SPDE driven by both Gaussian and non-Gaussian noise is much more difficult. This is due to the extra term appearing in the Burkholder-Davis-Gundy inequality \eqref{Davis4} for compensated Poisson measure, which  makes estimates not easily optimal for $p>2$.  Observe that this extra term does not appear  in the Burkholder-Davis-Gundy inequality \eqref{Davis3} for Gaussian noise. In fact, a standard direction to extend the result in \cite{Yang} to SPDEs with additive noises to schemes achieving convergence order greater than $1/2$ could be to follow the work in \cite{Xiaojie2,Xiaojie3,AntonioRev1,AntonioRev2} by applying the Taylor expansion of order $2$ to the term $III_{11}$ \footnote{ In equation \eqref{lun3}}. But this requires the time regularity of the mild solution in the space $L^p(\Omega,H)$ for $p\in\{2,4\}$. Note that the time regularity requires the Burkholder-Davis-Gundy inequality. Note that for $p>2$, due the extra term in \eqref{Davis4},  the time regularity in $L^p(\Omega, H)$ has an optimal order $1/p$, instead of $1/2$ for Gaussian noise, and this  leads to convergence order $1/2$. Here, to achieve convergence order greater than $1/2$ for Gaussian and non Gaussian noises, we follow a different approach and only apply Taylor expansion of order $1$.  Note  that even with only Gaussian noise, such result has not yet  been obtained so far in the literature, to be best of our knowledge.
Note that this work can  be extended to the case of time-fractional stochastic partial differential equation driven  by both Gaussian  and non-Gaussian noises, by following  \cite{Zou1,Zou2}. 

 The rest of the paper is organised as follows.  \secref{wellposed} deals with the mathematical setting and the well posedness problem. In \secref{optimalregularity}, we provide optimal regularity results of the mild solution. In \secref{spacediscretization}, we perform the space discretisation of the problem with the finite element method and its error analysis. \secref{fulldiscretization1} deals with the full discretisation with exponential Euler method  and  semi implicit  scheme and their strong convergence 
analyses. 
In \secref{numericalexperiments}, we present numerical experiments to illustrate our theoretical results.

 \section{Mathematical setting,  main assumptions and well posedness problem}
 \label{wellposed}
 For a Banach space $U$ and any $p\geq 2$, we denote by $L^p(\Omega, U)$ the Banach space of all equivalence classes of $L^p$-integrable $U$-valued random variables. Let $L(U,H)$ be 
 the space of bounded linear mappings from $U$ to $H$ endowed with the usual  operator norm $\Vert .\Vert_{L(U,H)}$. 
 By  $\mathcal{L}_2(U,H):=HS(U,H)$,
 we  denote the space of Hilbert-Schmidt operators from $U$ to $H$ equipped with the norm 
 \begin{eqnarray}
 \Vert l\Vert^2_{\mathcal{L}_2(U,H)}:=\sum\limits_{i=1}^{\infty}\Vert l\psi_i\Vert^2, \quad  l\in \mathcal{L}_2(U,H),
 \end{eqnarray}
  where $(\psi_i)_{i=1}^{\infty}$ is an orthonormal basis of $U$. Note that this definition is independent of the orthonormal basis of $U$.  
For simplicity, we use the notations $L(U,U)=:L(U)$. and $\mathcal{L}_2(U,U)=:\mathcal{L}_2(U)$. 
 For all $l\in L(U,H)$ and $l_1\in\mathcal{L}_2(U)$ we have $ll_1\in\mathcal{L}_2(U,H)$ and 
\begin{eqnarray}
\label{chow1}
\Vert ll_1\Vert_{\mathcal{L}_2(U,H)}\leq \Vert l\Vert_{L(U,H)}\Vert l_1\Vert_{\mathcal{L}_2(U)}.
\end{eqnarray} 
 The space of Hilbert-Schmidt operators from  $Q^{1/2}(H)$ to $H$ is denoted by $L^0_2:=\mathcal{L}_2(Q^{1/2}(H),H)=HS(Q^{1/2}(H),H)$. As usual, $L^0_2$ is equipped with the norm
 \begin{eqnarray}
 \Vert l\Vert_{L^0_2} :=\Vert lQ^{1/2}\Vert_{HS}=\left(\sum\limits_{i=1}^{\infty}\Vert lQ^{1/2}e_i\Vert^2\right)^{1/2}, \quad  l\in L^0_2,
 \end{eqnarray}
where $(e_i)_{i=1}^{\infty}$ is an orthonormal basis  of $H$.
This definition is independent of the orthonormal basis of $H$.
 Let $L^p_{\mathcal{F}, \nu}(\chi\times[0,T];H)$, $p\geq 2$ be the space of all mappings $\theta :\chi\times[0,T]\times\Omega\longrightarrow H$ such that $\theta$ is jointly measurable and $\mathcal{F}_t$-adapted  for all $z\in\chi$, $0\leq s\leq t\leq T$  satisfying
\begin{eqnarray}
\mathbb{E}\left[\int_0^t\int_{\chi}\Vert \theta(z,s)\Vert^p\nu(dz)ds\right]<\infty.
\end{eqnarray}
\begin{lemma}\textbf{[Burkh\"{o}lder-Davis-Gundy inequalities]}
\label{Davis2}
Let $p\in[2,\infty)$ and $t\in[0, T]$
\begin{itemize}
\item[(i)] For any $L^0_2$-valued predictable process $\phi(s)$, $s\in[0, t]$ the following estimate holds
\begin{eqnarray}
\label{Davis3}
\mathbb{E}\left[\left\Vert\int_0^t\phi(s)dW(s)\right\Vert^p\right]\leq C(p)\mathbb{E}\left[\left(\int_0^t\Vert \phi(s)\Vert^2_{L^0_2}ds\right)^{p/2}\right].
\end{eqnarray}
\item[(ii)] For any $\theta\in L^p_{\mathcal{F}, \nu}(\chi\times[0, t]; H)$, $t\in[0, T]$ the following estimate holds
\begin{eqnarray}
\label{Davis4}
\mathbb{E}\left[\left\Vert\int_0^t\int_{\chi}\theta(z,s)\widetilde{N}(dz, ds)\right\Vert^p\right]&\leq&C(p)\mathbb{E}\left[\left(\int_0^t\int_{\chi}\Vert \theta(z,s)\Vert^2\nu(dz)ds\right)^{p/2}\right]\nonumber\\
&+&C(p)\mathbb{E}\left[\int_0^t\int_{\chi}\Vert \theta(z,s)\Vert^p\nu(dz)ds\right].
\end{eqnarray}
\end{itemize}
\end{lemma}
\begin{proof}
 The proof of \eqref{Davis3} can be found in \cite[Theorem 4.36]{Prato} and the proof of \eqref{Davis4} can be found in  \cite{Hausenblas} or  \cite[Theorem 3.2]{Gundy1}. Note that in the case $p=2$, the constants $C(p)$ in \eqref{Davis3} and \eqref{Davis4} are $1$, the inequality \eqref{Davis3} becomes the equality, the second term of \eqref{Davis4} does not exists any more and the equality in \eqref{Davis4} is replaced by the equality. In that case, they are known as the It\^{o} isometry property, see e.g. \cite[(4.30)]{Prato} and \cite[(3.56)]{Rudiger}. 
\end{proof}

In order to ensure the existence of the unique solution of \eqref{model}, we make the following assumptions.
\begin{Assumption}\textbf{[Initial value $X_0$]}
\label{assumption1}
We assume the initial data $X_0$ to be $\mathcal{F}_0$-measurable and  $X_0\in L^2\left(\Omega, \mathcal{D}\left((-A)^{\beta/2}\right)\right)$ and $0\leq\beta\leq 2$.
\end{Assumption}
\begin{Assumption}\textbf{[Nonlinear term $F$]}
\label{assumption2} 
We assume  the nonlinear mapping  $F: H\longrightarrow H$ to be  measurable and Lipschitz continuous, i.e. there exists a constant  $L>0$ such that 
\begin{eqnarray}
\label{reviewinequal1}
\Vert F(0)\Vert\leq L,\quad \Vert F(u)-F(v)\Vert \leq L\Vert u-v\Vert, \quad u,v\in H.
\end{eqnarray}
\end{Assumption}
 \begin{Assumption}\textbf{[Diffusion term]} 
 \label{assumption3}
  We assume that the operator  $B : H \longrightarrow L^0_2$ satisfies the global Lipschitz condition, i.e. there exists a positive constant $C$ such that 
 \begin{eqnarray}
 \Vert B(0)\Vert_{L^0_2}\leq L,\quad \Vert B(u)-B(v)\Vert_{L_2^0}\leq C\Vert u-v\Vert, \quad u,v\in H.
 \end{eqnarray}
 \end{Assumption}
 \begin{Assumption}\textbf{[Jump function]}
 \label{assumption4}
 The function $G: \chi\times H\longrightarrow H$, usually called jump coefficient or jump function is assumed to be measurable and satisfies
 \begin{eqnarray}
 \int_{\chi}\Vert G(z,0)\Vert^2\nu(dz)\leq L,\quad \int_{\chi}\Vert G(z, u)-G(z, v)\Vert^2\nu(dz)\leq L\Vert u-v\Vert^2,\quad u, v\in H.
 \end{eqnarray}
\end{Assumption}
The following theorem provides the well posedness of the problem \eqref{model}.
\begin{theorem}\cite{Alberverio}
\label{existence}
Let Assumptions \ref{assumption1}-\ref{assumption4} be fulfilled. Then the SPDE \eqref{model} has a unique mild solution $X(t)$ satisfying
\begin{itemize}
\item[(i)] $X(t)$ is $\mathcal{F}_t$-adapted on the filtration $(\Omega, \mathcal{F}, (\mathcal{F})_{t\geq 0}, \mathbb{P})$,
\item[(ii)] $\{X(t), t\in[0, T]\}$ is measurable and $\mathbb{E}\left[\displaystyle\int_0^T\Vert X(t)\Vert^2dt\right]<\infty$,
\item[(iii)] For all $t\in[0, T]$, the following equation holds almost surely
\begin{eqnarray}
\label{mild1}
X(t)&=&S(t)X_0+\int_0^tS(t-s)F\left(X(s)\right)ds+\int_0^tS(t-s)B\left(X(s)\right)dW(s)\nonumber\\
&+&\int_0^t\int_{\chi}S(t-s)G\left(z,X(s))\right)\widetilde{N}(dz,ds).
\end{eqnarray}
\end{itemize}
Moreover, it holds that $\sup\limits_{t\in[0, T]}\Vert X(t)\Vert_{L^2(\Omega, H)}<\infty$.
\end{theorem}

\begin{corollary}
\label{coro1}
As a consequence of \thmref{existence} and Assumptions \ref{assumption1}-\ref{assumption4}, the following  holds
\begin{eqnarray}
\sup_{0\leq t\leq T}\Vert F(X(t))\Vert_{L^2(\Omega, H)}+\sup_{0\leq t\leq T}\Vert B(X(t))\Vert_{L^2(\Omega, H)}+\int_{\chi}\sup_{0\leq t\leq T}\Vert G(z, X(t))\Vert_{L^2(\Omega, H)}\nu(dz)<\infty.
\end{eqnarray}
\end{corollary}
\section{Optimal regularity of the mild solution} 
\label{optimalregularity}
We discuss the space and  regularities of the mild solution $X(t)$ of \eqref{model} given by \eqref{mild1} in this section. In the rest of this paper, to simplify the presentation, we assume the SPDE \eqref{model} to be of second order of the following type
\begin{eqnarray}
\label{model1}
dX(t, x)&=&\left[\nabla \cdot \left(\mathbf{D}\nabla X(t,x)\right)-\mathbf{q}\cdot\nabla X(t,x)+f(x,X(t,x))\right]dt+b(x, X(t,x))dW(t,x)\nonumber\\
&+&\int_{\chi}g(z,x,X(t,x))\widetilde{N}(dz,dt),\quad x\in\Lambda,\quad t\in[0,T],
\end{eqnarray}
where $f:\Lambda\times \mathbb{R}\longrightarrow\mathbb{R}$ is globally Lipschitz continuous, $b:\Lambda\times\mathbb{R}\longrightarrow\mathbb{R}$ is continuously differentiable with globally bounded derivatives and $g:\chi\times\Lambda\times\mathbb{R}\longrightarrow\mathbb{R}$ is globally Lipschitz continuous. 
In the abstract framework \eqref{model}, the linear operator $A$ takes the form
\begin{eqnarray}
\label{operator}
Au=\sum_{i,j=1}^{d}\dfrac{\partial}{\partial x_i}\left(D_{ij}(x)\dfrac{\partial u}{\partial x_j}\right)-\sum_{i=1}^dq_i(x)\dfrac{\partial u}{\partial x_i},\quad
\mathbf{D}=\left(D_{i,j} \right)_{1\leq i,j \leq d},\,\,\,\,\,\,\, \mathbf{q}=\left( q_i \right)_{1 \leq i \leq d}.
\end{eqnarray}
where $D_{ij}\in L^{\infty}(\Lambda)$, $q_i\in L^{\infty}(\Lambda)$. We assume that there exists a positive constant $c_1>0$ such that 
\begin{eqnarray*}
\sum_{i,j=1}^dD_{ij}(x)\xi_i\xi_j\geq c_1|\xi|^2, \quad  \xi\in \mathbb{R}^d,\quad x\in\overline{\Omega}.
\end{eqnarray*}
The functions $F : H\longrightarrow H$, $B : H\longrightarrow L^0_2$ and $G:\chi\times H\longrightarrow H$ are defined by 
\begin{eqnarray}
\label{nemystskii}
(F(v))(x)=f(x,v(x)), \quad (B(v)u)(x)=b(x,v(x)).u(x) \quad \text{and}\quad G(z,v)(x)=g(z, x, v(x)),
\end{eqnarray}
for all $x\in \Lambda$, $v\in H$, $u\in Q^{1/2}(H)$ and $z\in\chi$.
As in \cite{Antonio1,Suzuki}, we introduce two spaces $\mathbb{H}$ and $V$, such that $\mathbb{H}\subset V$; the two spaces depend on the boundary 
conditions of $\Lambda$ and the domain of the operator $A$. For  Dirichlet (or first-type) boundary conditions we take 
\begin{eqnarray*}
V=\mathbb{H}=H^1_0(\Lambda)=\{v\in H^1(\Lambda) : v=0\quad \text{on}\quad \partial \Lambda\}.
\end{eqnarray*}
For Robin (third-type) boundary condition and  Neumann (second-type) boundary condition, which is a special case of Robin boundary condition, we take $V=H^1(\Lambda)$
\begin{eqnarray*}
\mathbb{H}=\{v\in H^2(\Lambda) : \partial v/\partial \mathtt{v}_{ A}+\alpha_0v=0,\quad \text{on}\quad \partial \Lambda\}, \quad \alpha_0\in\mathbb{R},
\end{eqnarray*}
where $\partial v/\partial \mathtt{v}_{ A}$ is the normal derivative of $v$ and $\mathtt{v}_{ A}$ is the exterior pointing normal $n=(n_i)$ to the boundary of $A$ given by
\begin{eqnarray*}
\partial v/\partial\mathtt{v}_{A}=\sum_{i,j=1}^dn_i(x)D_{ij}(x)\dfrac{\partial v}{\partial x_j},\,\,\qquad x \in \partial \Lambda.
\end{eqnarray*}
As in \cite{Suzuki, Antonio1}, one can easily check that $A$ generates an analytic semigroup $S(t)=e^{At}$ on $L^2(\Lambda)$ and the fractional powers of $-A$ are well defined. 

The following lemma will be of great interest in the next sections. Although its proof for $\rho\in[0,1)$ can be easily obtained from the smoothing properties of the semigroup \cite[Proposition 2.2]{Antonio1} or \cite{Henry}, the border case $\rho=1$ is critical and useful to achieve optimal convergence rates. The proof in the case of self-adjoint operator can be found in \cite[Lemma 3.2 (iii)]{Stig1} and the case of not necessarily self adjoint case can be found in \cite[Lemma 2.1]{Antjd1}.
 \begin{lemma}
 \label{sharpestimates}
 Let  $\rho\in[0, 1]$, then the following sharp integral estimates hold
 \begin{eqnarray}
 \label{sharp1}
 \int_{\tau_1}^{\tau_2}\Vert (-A)^{\rho/2}S(\tau_2-s)\Vert_{L(H)}^2 ds\leq C(\tau_2-\tau_1)^{1-\rho},\quad 0\leq \tau_1\leq \tau_2\leq T,\\
 \label{sharp2}
 \left\Vert  \int_{\tau_1}^{\tau_2}(-A)^{\rho}S(\tau_2-s) ds\right\Vert_{L(H)}\leq C(\tau_2-\tau_1)^{1-\rho},\quad 0\leq \tau_1\leq \tau_2\leq T.
 \end{eqnarray}
 Moreover, \eqref{sharp1} and \eqref{sharp2} hold if $A$ and $S$ are replaced by their discrete versions $A_h$ and $S_h$ respectively; defined in \secref{spacediscretization}.
 \end{lemma}
 We will also be interested on SPDE \eqref{model} driven by additive Wiener and Poisson measure noises, that is the SPDE \eqref{model} takes the following form
 \begin{eqnarray}
 \label{model3}
 dX(t)=[AX(t)+F(X(t))]dt+dW(t)+\int_{\chi}\psi(z)\widetilde{N}(dz,dt),\quad X(0)=X_0,\quad z_0\in\chi.
  \end{eqnarray}
 When dealing with additive Wiener and Poisson measure noises, we make the following assumption.
 \begin{Assumption}
 \label{assumption5}
 The covariance operator $Q:H\longrightarrow H$ and the jump coefficient $\psi :\chi\longrightarrow H$  satisfy the following estimates
 \begin{eqnarray}
 \label{wiener}
 \Vert (-A)^{\frac{\beta-1}{2}}Q^{\frac{1}{2}}\Vert_{\mathcal{L}_2(H)}<\infty,\quad
\int_{\chi} \Vert (-A)^{\frac{\beta-1}{2}}\psi(z)\Vert^2\nu(dz)<\infty.
 \end{eqnarray}
 where $\beta$ is given in \assref{assumption1}.
 \end{Assumption}
Let us now investigate the space regularity of the mild solution.
 \begin{proposition}
 \label{spaceregular}
 Let Assumptions \ref{assumption1}-\ref{assumption4} be fulfilled. Then the following space regularity holds
 \begin{eqnarray}
 \label{space0}
 \Vert (-A)^{\beta/2}X(t)\Vert_{L^2(\Omega, H)}\leq C\left(1+\Vert (-A)^{\beta/2}X_0\Vert_{L^2(\Omega, H)} \right),\quad 0\leq t\leq T,\quad 0\leq \beta\leq 1.
 \end{eqnarray}
 In the case of additive Wiener  and Poisson measure noises, if Assumptions \ref{assumption1}, \ref{assumption2} and \ref{assumption5}  are fulfilled, then the following estimate holds
 \begin{eqnarray}
 \label{space00}
 \Vert (-A)^{\beta/2}X(t)\Vert_{L^2(\Omega, H)}\leq C\left(1+\Vert (-A)^{\beta/2}X_0\Vert_{L^2(\Omega, H)} \right),\quad 0\leq t\leq T,\quad 0\leq \beta\leq 2.
 \end{eqnarray}
 Moreover, \eqref{space0} and \eqref{space00} hold if $A$ and $X$ are replaced respectively  by their discrete versions $A_h$ and $X^h$ defined in \secref{spacediscretization}.
 \end{proposition}
 \begin{proof}
 The proof of \eqref{space0} for $0\leq\beta<1$ can be found in \cite[Theorem 3.1]{Yang}. The border case $\beta=1$, useful in numerical analysis as it allows to achieve optimal convergence order was not treated. This can be done by following some steps of \cite[Theorem 3.1]{Yang}, but  further by using \lemref{sharpestimates}. The proof of \eqref{space00} follows also the same steps. 
\end{proof}

The following theorem provides a temporal regularity of the solution process of \eqref{model}. 
\begin{theorem}
\label{timeregular}
Suppose that Assumptions \ref{assumption1}-\ref{assumption4} are fulfilled. Then the following estimate holds
\begin{eqnarray}
\label{time0}
\Vert X(t_2)-X(t_1)\Vert_{L^2(\Omega, H)}\leq C(t_2-t_1)^{\min(\beta/2, 1/2)},\quad 0\leq t_1\leq t_2\leq T.
\end{eqnarray}
In the case of additive Wiener  and Poisson noises, if Assumptions \ref{assumption1}, \ref{assumption2} and \ref{assumption5} are fulfilled, then \eqref{time0} holds.

Moreover, \eqref{time0}  holds when $A$ and $X$ are replaced by their discrete versions $A_h$ and $X^h$ respectively; defined in \secref{spacediscretization}.
\end{theorem}
\begin{proof}
The proof of \eqref{time0} follows the same lines as that of \cite[Theorem 3.2]{Yang}. Note that the optimal regularity rate $1/2$ was not achieved in \cite{Yang}. We achieve that optimal regularity estimate here with the help of \lemref{sharpestimates}, which leads finally to the optimal convergence rate of the numerical schemes.  To prove \eqref{time0} in the case of additive noises, it follows from the mild solution \eqref{mild1} that
\begin{eqnarray}
\label{time1}
X(t_2)-X(t_1)
&=&S(t_1)\left (S(t_2-t_1)-\mathbf{I}\right)X_0+\int_0^{t_1}\left(S(t_2-t_1)-\mathbf{I}\right)S(t_1-s)F(X(s))ds\nonumber\\
&+&\int_{t_1}^{t_2}S(t_2-s)F(X(s))ds+\int_{t_1}^{t_2}S(t_2-s)dW(s)\nonumber\\
&+&\int_0^{t_1}(S(t_2-t_1)-\mathbf{I})S(t_1-s)dW(s)\nonumber\\
&+&\int_0^{t_1}\int_{\chi}(S(t_2-t_1)-\mathbf{I})S(t_1-s)\psi(z)\widetilde{N}(dz,ds)\nonumber\\
&+&\int_{t_1}^{t_2}\int_{\chi}S(t_2-s)\psi(z)\widetilde{N}(dz,ds)=:\sum_{i=0}^6I_i.
\end{eqnarray}
Taking the $L^p$ norm in both sides of \eqref{time1} and using triangle inequality yields
\begin{eqnarray}
\label{time2}
\left\Vert X(t_2)-X(t_1)\right\Vert_{L^2(\Omega, H)}\leq \sum_{i=0}^6\Vert I_i\Vert_{L^2(\Omega, H)}.
\end{eqnarray}
Using \assref{assumption1} and the smoothing property of the semigroup in \cite[Proposition 2.2]{Antonio1} or \cite{Henry} yields
\begin{eqnarray}
\label{time3}
\Vert I_0\Vert_{L^2(\Omega, H)}\leq\left\Vert S(t_1)\right\Vert_{L(H)}\left\Vert(S(t_2-t_1)-\mathbf{I})(-A)^{-\beta/2}\right\Vert_{L(H)}\Vert (-A)^{\beta/2}X_0\Vert_{L^2(\Omega, H)}\leq C(t_2-t_1)^{\beta/2}.
\end{eqnarray}
Using again the smoothing properties of the semigroup, \coref{coro1} and \lemref{sharpestimates} yields
\begin{eqnarray}
\label{time4}
\Vert I_1\Vert_{L^2(\Omega, H)}&\leq& \int_0^{t_1}\left\Vert (-A)^{-1/2}(S(t_2-t_1)-\mathbf{I})\right\Vert_{L(H)}\Vert (-A)^{1/2}S(t_1-s)\Vert_{L(H)}\Vert F(X(s))\Vert_{L^2(\Omega, H)}ds\nonumber\\
&\leq& C(t_2-t_1)^{1/2}\int_0^{t_1}\left\Vert (-A)^{1/2}S(t_1-s)\right\Vert_{L(H)}ds\leq C(t_2-t_1)^{1/2}.
\end{eqnarray}
Along the same lines as that for $I_1$, it holds that
\begin{eqnarray}
\label{time5}
\Vert I_2\Vert_{L^2(\Omega, H)}\leq C(t_2-t_1).
\end{eqnarray}
Let us move to the estimate of $I_5$.
Using the It\^{o} isometry  (Burkholder-Davis-Gundy inequality with $p=2$)   yields
\begin{eqnarray}
\label{time8}
\Vert I_5\Vert^2_{L^2(\Omega, H)}&=& \mathbb{E}\left[\int_0^{t_1}\int_{\chi}\left\Vert \left(S(t_2-t_1)-\mathbf{I}\right)S(t_1-s)\psi(z)\right\Vert^2\nu(dz)ds\right].
\end{eqnarray}
In view of the smoothing properties of the semigroup, \assref{assumption5} and   \lemref{sharpestimates}, it holds that
\begin{eqnarray}
\label{time9}
\Vert I_5\Vert^2_{L^2(\Omega, H)}&\leq& \int_0^{t_1}\int_{\chi}\Vert(-A)^{\frac{-\beta}{2}}(S(t_2-t_1)-\mathbf{I})\Vert^2_{L(H)}\Vert (-A)^{\frac{1}{2}} S(t_1-s)\Vert^2_{L(H)}\nonumber\\
&&\times\Vert(-A)^{\frac{\beta-1}{2}}\psi(z)\Vert^2\nu(dz)ds\nonumber\\
&\leq& C(t_2-t_1)^{\beta}\int_0^{t_1}\Vert (-A)^{\frac{1}{2}} S(t_1-s)\Vert^2_{L(H)}ds\int_{\chi}\Vert(-A)^{\frac{\beta-1}{2}}\psi(z)\Vert^2\nu(dz)\nonumber\\
&\leq& C(t_2-t_1)^{\beta}.
\end{eqnarray}
Using the It\^{o} isometry  (Burkholder-Davis-Gundy inequality with $p=2$), in view of the smoothing properties of the semigroup, \assref{assumption5} and   \lemref{sharpestimates}, it holds 
\begin{eqnarray}
\label{time10}
\Vert I_6\Vert^2_{L^2(\Omega, H)}&=& \mathbb{E}\left[\int_{t_1}^{t_2}\int_{\chi}\left\Vert S(t_2-s)\psi(z)\right\Vert^2\nu(dz)ds\right].
\nonumber\\
&\leq& \int^{t_2}_{t_1}\int_{\chi}\Vert (-A)^{\frac{1-\beta}{2}} S(t_2-s)\Vert^2_{L(H)}\Vert(-A)^{\frac{\beta-1}{2}}\psi(z)\Vert^2\nu(dz)ds\nonumber\\
&\leq& \int^{t_2}_{t_1}(t_2-s)^{\min(0,\beta-1)}ds\int_{\chi}\Vert(-A)^{\frac{\beta-1}{2}}\psi(z)\Vert^2\nu(dz)\nonumber\\
&\leq& C(t_2-t_1)^{\min(1,\beta)}.
\end{eqnarray}
Along the same lines as that of $\Vert I_5\Vert^2_{L^2(\Omega,H)}$ and $\Vert I_6\Vert^2_{L^2(\Omega,H)}$, we obtain the following estimates
\begin{eqnarray}
\label{time7}
\Vert I_3\Vert^2_{L^2(\Omega, H)}\leq C\Delta t^{\min(\beta,1)},\quad \Vert I_4\Vert^2_{L^2(\Omega,H)}\leq C\Delta t^{\beta}.
\end{eqnarray}
Substituting \eqref{time10}, \eqref{time9}, \eqref{time7}, \eqref{time5}, \eqref{time4} and \eqref{time3} in \eqref{time2} yields
\begin{eqnarray}
\left\Vert X(t_2)-X(t_1)\right\Vert_{L^2(\Omega, H)}\leq C(t_2-t_1)^{\min(\beta, 1)/2}.
\end{eqnarray}
The proof of the lemma is completed. 
 \end{proof}
 \section{Space approximation and error estimates}
 \label{spacediscretization}
Let us now turn to the space discretization of our problem \eqref{model}.  We start by  splitting  the domain $\Lambda$ in finite triangles.
Let $\mathcal{T}_h$ be the triangulation with maximal length $h$ satisfying the usual regularity assumptions, and  $V_h \subset V$ be the space of continuous functions that are 
piecewise linear over the triangulation $\mathcal{T}_h$. 
We consider the projection $P_h$ from $H=L^2(\Lambda)$ to $V_h$ defined by
\begin{eqnarray}
\label{projection}
\langle P_hu,\chi\rangle_H=\langle u,\chi\rangle_H, \quad \forall \chi\in V_h,\quad u\in H.
\end{eqnarray}
The discrete operator $A_h : V_h\longrightarrow V_h$ is defined by 
\begin{eqnarray}
\label{discreteoperator}
\langle A_h\phi,\chi\rangle_H=\langle A\phi,\chi\rangle_H=a(\phi,\chi),\quad \forall \phi,\chi\in V_h,
\end{eqnarray}
Like $-A$, $-A_h$ is also a generator of an analytic semigroup $S_h(t) : =e^{-tA_h}$.
The semi-discrete problem associated to \eqref{model} consists of finding $X^h(t)=X^h(.,t)\in V_h$ such that $X^h_0=P_hX_0$ and 
\begin{eqnarray}
\label{model4}
dX^h(t)=[A_hX^h(t)+P_hF(X^h(t))]dt+P_hB(X^h(t))dW(t)+\int_{\chi}P_hG(z,X^h(t))\widetilde{N}(dz,dt).
\end{eqnarray}
Note that obviously $A_h$, $P_hF$, $P_hB$ and $P_hG$ satisfy the same assumptions as $A$, $F$, $B$ and $G$ respectively. The mild solution of \eqref{model4} is given by
\begin{eqnarray}
\label{mild2}
X^h(t)&=&S_h(t)X^h_0+\int_0^tS_h(t-s)P_hF(X^h(s))ds+\int_0^tS_h(t-s)P_hB(X^h(s))dW(s)\nonumber\\
&+&\int_0^t\int_{\chi}S_h(t-s)P_hG(z, X^h(s))\widetilde{N}(dz,ds),\quad t\in[0,T].
\end{eqnarray}
To establish  our $L^2$ strong convergence result when dealing with multiplicative Wiener and Poisson measure noises, we will also need the following further assumption on the diffusion and jump functions when $\beta \in [1,2)$, which was also used in  \cite{Antonio1,Raphael,Stig1} to achieve optimal convergence order in space and optimal regularity results.
\begin{Assumption}
\label{assumption6}
We assume that there exists a positive constant $c>0$ such 
that $B\left(\mathcal{D}(-A)^{\gamma/2}\right)\subset HS\left(Q^{1/2}(H),\mathcal{D}(-A)^{\gamma/2}\right)$ and $G\left(z, \mathcal{D}((-A)^{\gamma/2})\right)\subset\mathcal{D}\left((-A)^{\gamma/2})\right)$, $z\in \chi$ and
\begin{eqnarray}
\Vert (-A)^{\gamma/2}B(v)\Vert_{L^0_2}\leq c(1+\Vert v\Vert_{\gamma}),\quad \Vert (-A)^{\gamma/2}G(z,v)\Vert\leq c(1+\Vert v\Vert_{\gamma}),\; v\in \mathcal{D}((-A)^{\gamma/2}),\; z\in \chi,
\end{eqnarray}
with $\gamma=\beta-1$, where $\beta$ is the parameter defined in Assumption \ref{assumption2}.
\end{Assumption}
\begin{theorem}
\label{spaceapproximation}
Let $X$ and $X^h$ be the mild solutions of \eqref{model} and \eqref{model4} respectively. Suppose that Assumptions \ref{assumption1}-\ref{assumption4} and \ref{assumption6} are fulfilled. The following space error estimate holds
\begin{eqnarray}
\label{espace1}
\Vert X(t)-X^h(t)\Vert_{L^2(\Omega, H)}\leq Ch^{\beta},\quad t\in[0, T].
\end{eqnarray}
When dealing with additive Wiener and Poisson measure noises, if Assumptions \ref{assumption1}, \ref{assumption2} and \ref{assumption5} are fulfilled, then the error estimate \eqref{espace1} holds.
\end{theorem}
\begin{proof} The proof follows the same lines as \cite[Theorem 5.1]{Yang}. Note that to achieve optimal convergence order, we use \lemref{spaceregular}, \propref{sharpestimates} and \cite[Lemma 3.2]{Antjd2}.
\end{proof}
\section{Fully discrete schemes}
\label{fulldiscretization1}
Applying the exponential Euler method \cite{Antonio1} to the semi-discrete problem \eqref{model4} yields the following fully discrete scheme
\begin{eqnarray}
\label{schema1}
X^h_{m+1}&=&S_h(\Delta t)X^h_m+\int_{t_m}^{t_{m+1}}S_h(t_{m+1}-s)P_hF(X^h_m)ds+\int_{t_m}^{t_{m+1}}S_h(\Delta t)P_hB(X^h_m)dW(s)\nonumber\\
&+& \int_{t_m}^{t_{m+1}}\int_{\chi}S_h(\Delta t)P_h G(z, X^h_m)\widetilde{N}(dz, ds),\quad X^h_0=P_hX_0.
\end{eqnarray}
    We also investigate the semi implicit scheme \cite{Yang} for not necessarily self adjoint operator 
    \begin{eqnarray}
    \label{vend1}
Y^h_{m+1}&=&S_{h, \Delta t}Y^h_m+\int_{t_m}^{t_{m+1}}S_{h, \Delta t}P_hF(Y^h_m)ds+\int_{t_m}^{t_{m+1}}S_{h, \Delta t}P_hB(Y^h_m)dW(s)\nonumber\\
&+&\int_{t_m}^{t_{m+1}}\int_{\chi}S_{h, \Delta t}P_hG(z, Y^h_m)\widetilde{N}(dz, ds),
\end{eqnarray}
where 
\begin{eqnarray}
\label{vend1a}
S_{h, \Delta t}:=(\mathbf{I}-\Delta tA_h)^{-1}\quad \text{and}\quad  Y^h_0=P_hX_0.
\end{eqnarray}

 In addition to the above assumptions, when dealing with additive noises, we need the following extra assumptions, useful to achieve convergence rate greater than $1/2$.
\begin{Assumption}
\label{assumption7} 
The drift function $F$ is assumed to be differentiable with derivatives satisfying
\begin{eqnarray}
\Vert F'(u)v\Vert\leq C\Vert v\Vert\quad \Vert (-A)^{-\frac{\eta}{2}}\left(F'(u)-F'(v)\right)\Vert_{L(H)}\leq C\Vert u-v\Vert,\quad 
u,v\in H, 
\end{eqnarray}
for some $\eta\in\left(\frac{3}{4}, 1\right)$. 
\end{Assumption}
\begin{lemma}
\label{distlemma}
Under Assumptions \ref{assumption5} and \ref{assumption7}, the following estimates hold
\begin{eqnarray}
\label{distes1}
\Vert (-A_h)^{\frac{\beta-1}{2}}P_hQ^{\frac{1}{2}}\Vert_{\mathcal{L}_2(H)}\leq C,\quad \int_{\chi}\Vert (-A_h)^{\frac{\beta-1}{2}}P_h\psi(z)\Vert^2\leq C\quad \Vert P_hF'(u)v\Vert\leq C\Vert v\Vert,\quad u, v\in H.
\end{eqnarray}
 The following estimate also holds
 \begin{eqnarray*}
 \Vert (-A_h)^{-\frac{\eta}{2}}P_h\left(F'(u)-F'(v)\right)\Vert_{L(H)}\leq C\Vert u-v\Vert,\quad 
u,v\in H.
 \end{eqnarray*}
\end{lemma}
\begin{proof}
The proof  can be found in \cite[(63) \& (65)]{Antonio2}. 
\end{proof}

The strong convergence results of the fully discrete schemes are formulated in the following theorems.
\begin{theorem}\textbf{[Multiplicative Gaussian and Poisson measure noises]}
\label{mainresult1}
Let Assumptions \ref{assumption1}-\ref{assumption4} be fulfilled.  Let $\xi^h_m$ be the numerical approximation defined in \eqref{schema1}  or \eqref{vend1} ($\xi^h_m=X^h_m$ for exponential scheme \eqref{schema1}  and  $\xi^h_m=Y^h_m$ 
for  semi implicit scheme \eqref{vend1}). We have the following estimates depending of the regularity of the initial solution $X_0$
\begin{itemize}
\item[(i)] If $0<\beta<1$ then the following error estimate holds
\begin{eqnarray}
\Vert X(t_m)-\xi^h_m\Vert_{L^2(\Omega, H)}\leq C\left(h^{\beta}+\Delta t^{\beta/2}\right).
\end{eqnarray}
\item[(ii)] If $\beta\in[1, 2]$ then the following error estimate holds
\begin{eqnarray}
\Vert X^h(t_m)-\xi^h_m\Vert_{L^2(\Omega, H)}\leq C\left(h^{1-\epsilon}+\Delta t^{1/2}\right).
\end{eqnarray}
\item[(iii)] If $\beta\in[1, 2]$ and moreover if \assref{assumption5} is fulfilled, then the following error estimate holds
\begin{eqnarray}
\Vert X(t_m)-\xi^h_m\Vert_{L^2(\Omega, H)}\leq C\left(h^{\beta}+\Delta t^{1/2}\right).
\end{eqnarray}
\end{itemize}
\end{theorem}
\begin{theorem}\textbf{[Additive Gaussian and Poisson measure noises]}
\label{mainresult2}
Let Assumptions \ref{assumption1}, \ref{assumption6} and \ref{assumption7} be fulfilled.  Let $X^h_m$ and $Y^h_m$  be the numerical approximations defined in \eqref{schema1}  and \eqref{vend1} respectively.
Depending of the regularity of the initial solution $X_0$, the following error estimate holds
\begin{itemize}
\item [(i)]  If $\beta\in[0, 2]$, then the following holds
\begin{eqnarray}
\Vert X(t_m)-X^h_m\Vert_{L^2(\Omega, H)}\leq C\left(h^{\beta}+\Delta t^{\beta/2}\right).
\end{eqnarray}
\item[(ii)] If $\beta\in[0, 2)$, then the following holds
\begin{eqnarray}
\Vert X(t_m)-Y^h_m\Vert_{L^2(\Omega, H)}\leq C\left(h^{\beta}+\Delta t^{\beta/2}\right).
\end{eqnarray}
\item[(iii)] If $\beta=2$, then the following holds
\begin{eqnarray}
\Vert X(t_m)-Y^h_m\Vert_{L^2(\Omega, H)}\leq C\left(h^2+\Delta t^{1-\epsilon}\right).
\end{eqnarray}
\end{itemize}
\end{theorem}

\begin{remark}
\label{remarkconditions}
Note that the convergence order in time greater than $1/2$ was achieved in the literature for additive Gaussian noise by incorporating additional assumption on the second derivative of $F$, see e.g. \cite[Assumption 2.2]{Xiaojie2}, \cite[Assumption 2.1]{Xiaojie3}, \cite[Assumption 2.2]{AntonioRev1} and \cite[Assumption 2.1]{AntonioRev2}. Note that such assumptions on the nonlinearity $F$ is more restrictive since in some realistic situations the second derivative of $F$ may not exist. In this paper, we provide less restrictive assumptions of $F$, relying only on  its first derivative. Note that our relaxed \assref{assumption7} allows a large class of nonlinear Nemytskii operators. \thmref{mainresult2} shows that under \assref{assumption7}, we still achieve convergence order greater than $1/2$.
\end{remark}

\begin{remark}
\label{remarkorder}
In  \thmref{mainresult2}, we achieve optimal convergence rates $1$ and $1-\epsilon$ in time  for exponential  integrator and semi implicit scheme respectively.
 This extends the well known results in the case of stochastic differential equations driven by additive Brownian motion and Poisson process to the case of SPDE driven by Wiener process and Poisson random measure. 
 In fact,  for such problems for SDEs, since the diffusion  and the jump functions are constants,  their derivatives vanish and the Milstein schemes read as the  Euler schemes  therefore have  convergence order $1$.
\end{remark}
The proofs of Theorems \ref{mainresult1} and \ref{mainresult2} are given in the following sections.
\subsection{Proof of \thmref{mainresult1} for the exponential scheme}
Using trinagle inequality yields 
\begin{eqnarray}
\Vert X^h(t_m)-X^h_m\Vert_{L^2(\Omega, H)}\leq \Vert X(t_m)-X^h(t_m)\Vert_{L^2(\Omega, H)}+\Vert X^h(t_m)-X^h_m\Vert_{L^2(\Omega, H)}.
\end{eqnarray}
The space error is estimated in \thmref{spaceapproximation}. It remains to estimate the time error. Let us recall that the mild solution \eqref{mild2} can be written as follows
\begin{eqnarray}
\label{mild3}
X^h(t_m)&=&S_h(\Delta t)X^h(t_{m-1})+\int_{t_{m-1}}^{t_m}S_h(t_m-s)P_hF(X^h(s))ds\\
&+&\int_{t_{m-1}}^{t_m}S_h(t_m-s)P_hB(X^h(s))dW(s)+\int_{t_{m-1}}^{t_m}\int_{\chi}S_h(t_m-s)P_hG(z, X^h(s))\widetilde{N}(dz, ds).\nonumber
\end{eqnarray}
 Iterating the mild solution \eqref{mild3} yields
\begin{eqnarray}
\label{main1}
X^h(t_m)&=& S_h(t_m)P_hX_0+\sum_{k=0}^{m-1}\int_{t_k}^{t_{k+1}}S_h(t_{m}-s)P_hF(X^h(s))ds\nonumber\\
&+&\sum_{k=0}^{m-1}\int_{t_k}^{t_{k+1}}S_h(t_{m}-s)P_hB(X^h(s))dW(s)\nonumber\\
&+&\sum_{k=0}^{m-1}\int_{t_k}^{t_{k+1}}\int_{\chi}S_h(t_{m}-s)P_h G(z, X^h(s))\widetilde{N}(dz, ds).
\end{eqnarray}
Iterating the numerical solution \eqref{schema1} yields 
\begin{eqnarray}
\label{main2}
X^h_m&=& S_h(t_m)P_hX_0+\sum_{k=0}^{m-1}\int_{t_k}^{t_{k+1}}S_h(t_{m}-s)P_hF(X^h_k)ds+\sum_{k=0}^{m-1}\int_{t_k}^{t_{k+1}}S_h(t_m-t_k)P_hB(X^h_k)dW(s)\nonumber\\
&+&\sum_{k=0}^{m-1}\int_{t_k}^{t_{k+1}}\int_{\chi}S_h(t_m-t_k)P_h  G(z, X^h_k)\widetilde{N}(dz, ds).
\end{eqnarray}
Subtracting \eqref{main2} from \eqref{main1} yields
\begin{eqnarray}
\label{main3}
X^h(t_m)-X^h_m&=&\sum_{k=0}^{m-1}\int_{t_k}^{t_{k+1}}S_h(t_m-s)[P_hF(X^h(s))-P_hF(X^h_k)]ds\nonumber\\
&+&\sum_{k=0}^{m-1}\int_{t_k}^{t_{k+1}}[S_h(t_{m}-s)P_hB(X^h(s))-S_h(t_m-t_k)P_hB(X^h_k)]dW(s)\nonumber\\
&+&\sum_{k=0}^{m-1}\int_{t_k}^{t_{k+1}}\int_{\chi}[S_h(t_{m}-s)P_hG(z,X^h(s))-S_h(t_m-t_k)P_h G(z,X^h_k)]\widetilde{N}(dz,ds)\nonumber\\
&=:&II_1+II_2+II_3.
\end{eqnarray}
Using triangle inequality yields
\begin{eqnarray}
\label{main4}
\Vert X^h(t_m)-X^h_m\Vert_{L^2(\Omega, H)}\leq \Vert II_1\Vert_{L^2(\Omega, H)}+\Vert II_2\Vert_{L^2(\Omega, H)}+\Vert II_3\Vert_{L^2(\Omega, H)}.
\end{eqnarray}
We can recast $II_1$ as follows
\begin{eqnarray}
\label{main5a}
II_1&=&\sum_{k=0}^{m-1}\int_{t_k}^{t_{k+1}}S_h(t_{m}-s)[P_hF(X^h(s))-P_hF(X^h(t_k))]ds\nonumber\\
&+&\sum_{k=0}^{m-1}\int_{t_k}^{t_{k+1}}S_h(t_{m}-s)[P_hF(X^h(t_k))-P_hF(X^h_k)]ds=:II_{11}+II_{12}.
\end{eqnarray}
Using the smoothing properties of the semigroup, \assref{assumption2} and \eqref{time0} yields
\begin{eqnarray}
\label{main5b}
\Vert II_{11}\Vert_{L^2(\Omega, H)}\leq C\sum_{k=0}^{m-1}\int_{t_k}^{t_{k+1}}\Vert X^h(s)-X^h(t_k)\Vert_{L^2(\Omega, H)}\leq C\Delta t^{\min(\beta, 1)/2}.
\end{eqnarray}
In view of the stability property of the semigroup and \assref{assumption2}, it holds that
\begin{eqnarray}
\label{main5c}
\Vert II_{12}\Vert_{L^2(\Omega, H)}\leq C\Delta t\sum_{k=0}^{m-1}\Vert X^h(t_k)-X^h_k\Vert_{L^2(\Omega, H)}.
\end{eqnarray}
Substituting \eqref{main5c} and \eqref{main5b} in \eqref{main5a} yields
\begin{eqnarray}
\label{main5}
\Vert II_1\Vert_{L^2(\Omega, H)}\leq C\Delta t^{\min(\beta, 1)/2}+C\Delta t\sum_{k=0}^{m-1}\Vert X^h(t_k)-X^h_k\Vert_{L^2(\Omega, H)}.
\end{eqnarray}
We will not give details of the estimate of $II_2$ as it is similar to that  of $II_3$.
Let us now  estimate the norm of $II_3$. We can recast $II_3$ in three terms as follows
\begin{eqnarray}
\label{main6a}
II_3&=&\sum_{k=0}^{m-1}\int_{t_k}^{t_{k+1}}\int_{\chi}\left[S_h(t_{m}-s)-S_h(t_{m}-t_{m-1})\right]P_hG(z, X^h(s))\widetilde{N}(dz, ds)\nonumber\\
&+&\sum_{k=0}^{m-1}\int_{t_k}^{t_{k+1}}\int_{\chi}S_h(t_m-t_k)[P_hG(z, X^h(s))-P_hG(z, X^h(t_k))]\widetilde{N}(dz, ds)\nonumber\\
&+& \sum_{k=0}^{m-1}\int_{t_k}^{t_{k+1}}\int_{\chi}S_h(t_m-t_k)[P_hG(z, X^h(t_k))-P_hG(z, X^h_k)]\widetilde{N}(dz, ds)\nonumber\\
&=:&II_{31}+II_{32}+II_{33}.
\end{eqnarray} 
Applying the It\^{o} isometry property \cite[(3.5.6)]{Rudiger}, using the smoothing properties of the semigroup, the $\sigma$ finiteness of the measure $\nu$, \assref{assumption4}, \lemref{sharpestimates} and \coref{coro1} yields
\begin{eqnarray}
\label{main6b}
\Vert II_{31}\Vert^2_{L^2(\Omega, H)}
&=&\sum_{k=0}^{m-1}\mathbb{E}\int_{t_k}^{t_{k+1}}\int_{\chi}\Vert S_h(t_m-s)(\mathbf{I}-S_h(s-t_k))P_hG(z, X^h(s))\Vert^2\nu(dz, ds)\nonumber\\
&\leq& \sum_{k=0}^{m-1}\int_{t_k}^{t_{k+1}}\Vert (-A_h)^{1/2}S_h(t_m-s)\Vert^2_{L(H)}\Vert (-A_h)^{-1/2}(\mathbf{I}-S_h(s-t_k))\Vert^2_{L(H)}\nonumber\\
&&\times\mathbb{E}\int_{\chi}\Vert P_hG(z, X^h(s))\Vert^2\nu(dz)ds\nonumber\\
&\leq& C\sum_{k=0}^{m-1}\int_{t_k}^{t_{k+1}}(s-t_k)\Vert(-A_h)^{1/2}S_h(t_m-s)\Vert^2_{L(H)}\left(1+\Vert X^h(s)\Vert^2_{L(H)}\right)ds\nonumber\\
&\leq& C\Delta t\sum_{k=0}^{m-1}\int_{t_k}^{t_{k+1}}\Vert(-A_h)^{1/2}S_h(t_m-s)\Vert^2_{L(H)}ds\nonumber\\
&=&C\Delta t\int_0^{t_m}\Vert(-A_h)^{1/2}S_h(t_m-s)\Vert^2_{L(H)}ds\leq C\Delta t.
\end{eqnarray}
Applying again the It\^{o} isometry property \cite[(3.5.6)]{Rudiger}, using  the smoothing properties of the semigroup, \assref{assumption4} and \eqref{time0} yields
\begin{eqnarray}
\label{main6c}
\Vert II_{32}\Vert^2_{L^2(\Omega, H)}&=&\sum_{k=0}^{m-1}\mathbb{E}\int_{t_k}^{t_{k+1}}\int_{\chi}\Vert S_h(t_m-t_k)[P_hG(z, X^h(s))-P_hG(z, X^h(t_k))]\Vert^2\nu(dz)ds\nonumber\\
&\leq& C\sum_{k=0}^{m-1}\int_{t_k}^{t_{k+1}}\mathbb{E}\Vert X^h(s)-X^h(t_k)\Vert^2ds\nonumber\\
&\leq& C\sum_{k=0}^{m-1}\int_{t_k}^{t_{k+1}}(s-t_k)^{\min(\beta, 1)}ds\leq C\Delta t^{\min(\beta,1)}.
\end{eqnarray}
Applying the It\^{o} isometry property,  using the smoothing properties of the semigroup and  \assref{assumption4}  yields
\begin{eqnarray}
\label{main6d}
\Vert II_{33}\Vert^2_{L^2(\Omega, H)}&=&\sum_{k=0}^{m-1}\mathbb{E}\int_{t_k}^{t_{k+1}}\int_{\chi}\Vert S_h(t_m-t_k)[P_hG(z, X^h(t_k))-P_hG(z, X^h_k)]\Vert^2\nu(dz)ds\\
&\leq& C\sum_{k=0}^{m-1}\int_{t_k}^{t_{k+1}}\mathbb{E}\Vert X^h(t_k)-X^h_k\Vert^2ds\leq C\Delta t\sum_{k=0}^{m-1}\Vert X^h(t_k)-X^h_k\Vert^2_{L^2(\Omega, H)}.\nonumber
\end{eqnarray}
Substituting \eqref{main6d}, \eqref{main6c} and \eqref{main6b} in \eqref{main6a} yields
\begin{eqnarray}
\label{main6}
\Vert II_3\Vert^2_{L^2(\Omega, H)}\leq C\Delta t^{\min(\beta, 1)}+C\Delta t\sum_{k=0}^{m-1}\Vert X^h(t_k)-X^h_k\Vert^2_{L^2(\Omega, H)}.
\end{eqnarray}
Using the similar procedure as  for $II_3$, we obtain the following estimate of the norm of $II_2$
\begin{eqnarray}
\label{main7}
\Vert II_2\Vert^2_{L^2(\Omega, H)}\leq C\Delta t^{\min(\beta, 1)}+C\Delta t\sum_{k=0}^{m-1}\Vert X^h(t_k)-X^h_k\Vert^2_{L^2(\Omega, H)}.
\end{eqnarray}
Substituting \eqref{main7}, \eqref{main6} and \eqref{main5} in \eqref{main2} yields
\begin{eqnarray}
\label{main8}
\Vert X^h(t_m)-X^h_m\Vert^2_{L^2(\Omega, H)}\leq C\Delta t^{\min(\beta, 1)}+C\Delta t\sum_{k=0}^{m-1}\Vert X^h(t_k)-X^h_k\Vert^2_{L^2(\Omega, H)}.
\end{eqnarray}
Applying the discrete Gronwall lemma to \eqref{main8} and taking the square root yields
\begin{eqnarray}
\label{main9}
\Vert X^h(t_m)-X^h_m\Vert_{L^2(\Omega, H)}\leq C\Delta t^{\min(\beta, 1)/2}.
\end{eqnarray}
Combining \eqref{main9} and \thmref{spaceapproximation} completes the proof of \thmref{mainresult1}.

\subsection{Proof of \thmref{mainresult2} for exponential scheme}
Since the space error is estimated in \thmref{spaceapproximation}, we only need to estimate the time error. Note that form \eqref{main4}, we have
\begin{eqnarray}
\label{lun1}
\Vert X^h(t_m)-X^h_m\Vert_{L^2(\Omega, H)}\leq \Vert III_1\Vert_{L^2(\Omega, H)}+\Vert III_2\Vert_{L^2(\Omega, H)}+\Vert III_3\Vert_{L^2(\Omega, H)},
\end{eqnarray}
where $III_1$ is exactly the same as $II_1$ in \eqref{main3}. The terms involving the noises $III_2$ and $III_3$ are given by
\begin{eqnarray}
\label{lun2}
III_2&=&\sum_{k=0}^{m-1}\int_{t_k}^{t_{k+1}}\left[S_h(t_m-s)-S_h(t_m-t_k)\right]P_hdW(s),\\
III_3&=&\sum_{k=0}^{m-1}\int_{t_k}^{t_{k+1}}\int_{\chi}\left(S_h(t_m-s)-S_h(t_m-t_k)\right)P_h\psi(z)\widetilde{N}(dz, ds).
\end{eqnarray}
Let us start with the $L^2$ norm estimate of $III_1$. We recall that from \eqref{main5a} we have
\begin{eqnarray}
\label{lun3}
\Vert III_1\Vert_{L^2(\Omega, H)}\leq \Vert III_{11}\Vert_{L^2(\Omega, H)}+\Vert III_{12}\Vert_{L^2(\Omega, H)}, 
\end{eqnarray}
where $III_{11}$ and $III_{12}$ are  the same as $II_{11}$ and $II_{12}$ respectively. Note that from \eqref{main5c} we have
\begin{eqnarray}
\label{lun4}
\Vert III_{12}\Vert_{L^2(\Omega, H)}\leq C\Delta t\sum_{k=0}^{m-1}\Vert X^h(t_k)-X^h_k\Vert_{L^2(\Omega, H)}.
\end{eqnarray}
The standard estimate of $III_{11}$ as in \eqref{main5b}  gives a convergence order $1/2$.
 To achieve higher order, we need to re-estimate $III_{11}$. We are interested in the case $\beta\in(1,2]$, since the case of $\beta\in[0, 1]$ goes along the same lines as in the case multiplicative noise.  Let us recall that
\begin{eqnarray}
\label{lun5}
III_{11}=\sum_{k=0}^{m-1}\int_{t_k}^{t_{k+1}}S_h(t_m-s)P_h[F(X^h(s))-F(X^h(t_k))]ds.
\end{eqnarray}
We use the Taylor expansion in Banach space as in \cite{Jentzen1} to obtain
\begin{eqnarray}
\label{Taylor1}
F(X^h(s))-F(X^h(t_k))=\left(\int_0^1F'\left(X^h(t_k))+\lambda\left(X^h(s)-X^h(t_k)\right)\right)d\lambda\right)\left(X^h(s)-X^h(t_k)\right).
\end{eqnarray}
Note that the mild representation $X^h(s)$ can be written as follows
\begin{eqnarray}
\label{Taylor1a}
X^h(s)&=&S_h(s-t_k)X^h(t_k)+\int_{t_k}^sS_h(s-r)P_hF(X^h(r))dr+\int_{t_k}^sS_h(s-r)P_hdW(r)\nonumber\\
&+&\int_{t_k}^s\int_{\chi}S_h(s-r)P_h\psi(z)\widetilde{N}(dz,dr),\quad t_k\leq s.
\end{eqnarray}
Substituting \eqref{Taylor1a} in \eqref{Taylor1}, it follows that
\begin{eqnarray}
\label{Taylor2}
&&F(X^h(s))-F(X^h(t_k))\nonumber\\
&=&I^h_{k,s}\left(S_h(s-t_k)-\mathbf{I}\right)X^h(t_k)+I^h_{k,s}\int_{t_k}^sS_h(s-r)P_hF(X^h(r))dr\nonumber\\
&+&I^h_{k,s}\int_{t_k}^sS_h(s-r)P_hdW(r)+I^h_{k,s}\int_{t_k}^s\int_{\chi}S_h(s-r)P_h\psi(z)\widetilde{N}(dz,dr),
\end{eqnarray}
where
\begin{eqnarray}
\label{Taylor2a}
I^h_{k,s}:=\int_0^1F'\left(X^h(t_k)+\lambda\left(X^h(s)-X^h(t_k)\right)\right)d\lambda.
\end{eqnarray}
Substituting \eqref{Taylor2} in \eqref{lun5} yields
\begin{eqnarray}
\label{Taylor3}
III_{11}&=&\sum_{k=0}^{m-1}\int_{t_k}^{t_{k+1}}S_h(t_m-s)I^h_{k,s}\left(S_h(s-t_k)-\mathbf{I}\right)X^h(t_k)ds\nonumber\\
&+&\sum_{k=0}^{m-1}\int_{t_k}^{t_{k+1}}S_h(t_m-s)I^h_{k,s}\int_{t_k}^sS_h(s-r)P_hF(X^h(r))drds\nonumber\\
&+&\sum_{k=0}^{m-1}\int_{t_k}^{t_{k+1}}S_h(t_m-s)I^h_{k,s}\int_{t_k}^sS_h(s-r)P_hdW(r)ds\nonumber\\
&+&\sum_{k=0}^{m-1}\int_{t_k}^{t_{k+1}}S_h(t_m-s)I^h_{k,s}\int_{t_k}^s\int_{\chi}S_h(s-r)P_h\psi(z)\widetilde{N}(dz,dr)ds\nonumber\\
&=:&III_{11}^{(1)}+III_{11}^{(2)}+III_{11}^{(3)}+III_{11}^{(4)}.
\end{eqnarray}
Taking the norm in both sides of \eqref{Taylor3} and using triangle inequality yields
\begin{eqnarray}
\label{Taylor4}
\Vert III_{11}\Vert_{L^2(\Omega, H)}\leq \Vert III_{11}^{(1)}\Vert_{L^2(\Omega,H)}+\Vert III_{11}^{(2)}\Vert_{L^2(\Omega, H)}+\Vert III_{11}^{(3)}\Vert_{L^2(\Omega,H)}+\Vert III_{11}^{(4)}\Vert_{L^2(\Omega,H)}.
\end{eqnarray}
Using triangle inequality, the smoothing properties of the semigroup, Lemmas \ref{distlemma} and \ref{spaceregular} yields
\begin{eqnarray}
\label{Taylor5}
\Vert III_{11}^{(1)}\Vert_{L^2(\Omega,H)}&\leq& \sum_{k=0}^{m-1}\int_{t_k}^{t_{k+1}}\Vert S_h(t_m-s)I^h_{k,s}\left(S_h(s-t_k)-\mathbf{I}\right)X^h(t_k)\Vert_{L^2(\Omega, H)}ds\nonumber\\
&\leq& \sum_{k=0}^{m-1}\int_{t_k}^{t_{k+1}}\Vert S_h(t_m-s)I^h_{k,S}\Vert_{L(H)}\left\Vert\left(S_h(s-t_k)-\mathbf{I}\right)(-A_h)^{-\beta/2}\right\Vert_{L(H)}\nonumber\\
&&\times\left\Vert(-A_h)^{\beta/2}X^h(t_k)\right\Vert_{L^2(\Omega,H)}ds\nonumber\\
&\leq& C\sum_{k=0}^{m-1}\int_{t_k}^{t_{k+1}}(s-t_k)^{\beta/2}ds\leq C\Delta t^{\beta/2}.
\end{eqnarray}
Using triangle inequality, the smoothing properties of the semigroup, \lemref{distlemma} and \coref{coro1}, it holds that
\begin{eqnarray}
\label{Taylor6}
\Vert III_{11}^{(2)}\Vert_{L^2(\Omega,H)}&\leq& \sum_{k=0}^{m-1}\int_{t_k}^{t_{k+1}}\int_{t_k}^s\Vert S_h(t_m-s)I^h_{k,s}S_h(s-r)P_hF(X^h(r))\Vert_{L^2(\Omega,H)}drds\nonumber\\
&\leq& C\sum_{k=0}^{m-1}\int_{t_k}^{t_{k+1}}(s-t_k)ds\leq C\Delta t.
\end{eqnarray}
We split $III_{11}^{(3)}$ in two terms as follows
\begin{eqnarray}
\label{Raphael1}
III_{11}^{(3)}&=&\sum_{k=0}^{m-1}\int_{t_k}^{t_{k+1}}S_h(t_m-s)I^h_{k,k}\int_{t_k}^sS_h(s-r)P_hdW(r)ds\nonumber\\
&+&\sum_{k=0}^{m-1}\int_{t_k}^{t_{k+1}}S_h(t_m-s)\left(I^h_{k,s}-I^h_{k,k}\right)\int_{t_k}^sS_h(s-r)P_hdW(r)ds\nonumber\\
&=:&III_{11}^{(31)}+III_{11}^{(32)}.
\end{eqnarray}
Using the fact that the expectation of the cross-product terms vanishes, using H\"{o}lder inequality, the It\^{o}-isometry  and \lemref{distlemma}, the smoothing properties of the semigroup, it follows that
\begin{eqnarray}
\label{Taylor7}
\Vert III_{11}^{(31)}\Vert_{L^2(\Omega,H)}^2&=&\mathbb{E}
\left[\left\Vert\sum_{k=0}^{m-1}\int_{t_k}^{t_{k+1}}\left(\int_{t_k}^sS_h(t_m-s)I^h_{k,k}S_h(s-r)P_hdW(r)\right)ds\right\Vert^2\right]\nonumber\\
&=&\sum_{k=0}^{m-1}\mathbb{E}
\left\Vert\int_{t_k}^{t_{k+1}}\left(\int_{t_k}^sS_h(t_m-s)I^h_{k,k}S_h(s-r)P_hdW(r)\right)ds\right\Vert^2\nonumber\\
&\leq& \Delta t\sum_{k=0}^{m-1}\int_{t_k}^{t_{k+1}}\mathbb{E}\left\Vert\int_{t_k}^sS_h(s-t_k)I^h_{k,k}S_h(s-r)P_hdW(r)\right\Vert^2ds\nonumber\\
&\leq& \Delta t\sum_{k=0}^{m-1}\int_{t_k}^{t_{k+1}}\int_{t_k}^s\left\Vert S_h(t_m-s)I^h_{k,k}S_h(s-r)P_hQ^{\frac{1}{2}}\right\Vert^2_{\mathcal{L}_2(H)}drds\nonumber\\
&\leq& \Delta t\sum_{k=0}^{m-1}\int_{t_k}^{t_{k+1}}\int_{t_k}^s\Vert S_h(t_m-s)I^h_{k,k}\Vert^2_{L(H)}\left\Vert S_h(s-r)(-A_h)^{\frac{1-\beta}{2}}\right\Vert^2_{L(H)}\nonumber\\
&&\times\left\Vert (-A_h)^{\frac{\beta-1}{2}}P_hQ^{\frac{1}{2}}\right\Vert^2_{\mathcal{L}_2(H)}drds\nonumber\\
&\leq& C\Delta t\sum_{k=0}^{m-1}\int_{t_k}^{t_{k+1}}\int_{t_k}^sdrds\leq C\Delta t^2.
\end{eqnarray}
Using triangle inequality and Cauchy-Schwartz inequality yields
 \begin{eqnarray*}
 \Vert III_{11}^{(32)}\Vert^2_{L^2(\Omega, H)}
 &\leq& m\sum_{k=0}^{m-1}\left\Vert\int_{t_{k}}^{t_{k+1}}S_h(t_m-s)P_h\left(I^h_{k,s}-I^h_{k,k}\right)\int_{t_{k}}^sS_h(s-r)P_hdW(r)ds\right\Vert^2_{L^2(\Omega, H)}\nonumber\\
 &\leq& m\Delta t\sum_{k=0}^{m-1}\int_{t_{k}}^{t_{k+1}}\left\Vert\int_{t_{k}}^sS_h(t_m-s)P_h\left(I^h_{k,s}-I^h_{k,k}\right) S_h(s-r)P_hdW(r)\right\Vert^2_{L^2(\Omega, H)}ds\nonumber\\
 &\leq& T \sum_{k=0}^{m-1}\int_{t_{k}}^{t_{k+1}}\mathbb{E}
 \left\Vert S_h(t_m-s)P_h\left(I^h_{k,s}-I^h_{k,k}\right)\int_{t_{k}}^s S_h(s-r)P_hdW(r)\right\Vert^2ds\nonumber\\
 &\leq& T \sum_{k=0}^{m-1}\int_{t_{k}}^{t_{k+1}}\left(\mathbb{E}
 \left\Vert S_h(t_m-s)P_h\left(I^h_{k,s}-I^h_{k,k}\right)\right\Vert^4_{L(H)}\right)^{\frac{1}{2}}\nonumber\\
 &\times&\left(\mathbb{E}\left\Vert\int_{t_{k}}^s  S_h(s-r)P_hdW(r)\right\Vert^4\right)^{\frac{1}{2}} ds.
 \end{eqnarray*}
Using the Burkh\"{o}lder-Davis-Gundy inequality (\cite[Lemma 5.1]{Raphael})
 \begin{eqnarray}
 \label{mach5}
 \Vert III_{11}^{(32)}\Vert^2_{L^2(\Omega, H)}
 &\leq& T \sum_{k=0}^{m-1}\int_{t_{k}}^{t_{k+1}}\left(\mathbb{E}
 \left\Vert S_h(t_m-s)P_h\left(I^h_{k,s}-I^h_{k,k}\right)\right\Vert^4_{L(H)}\right)^{\frac{1}{2}}\nonumber\\
 &\times&\int_{t_{k}}^s \mathbb{E}\left\Vert S_h(s-r)P_hQ^{\frac{1}{2}}\right\Vert^2_{\mathcal{L}_2(H)} dr ds.
 \end{eqnarray}
 Using \assref{assumption5} and the smoothing properties of the semigroup and the fact that $\beta\in(1, 2]$, it holds that
 \begin{eqnarray}
 \label{mach6a}
 \left\Vert S_h(s-r)P_hQ^{\frac{1}{2}}\right\Vert^2_{\mathcal{L}_2(H)}
 &\leq&\left\Vert S_h(s-r)(-A_h(r))^{\frac{1-\beta}{2}}(-A(r))^{\frac{\beta-1}{2}} P_hQ^{\frac{1}{2}}\right\Vert^2_{\mathcal{L}_2(H)}\nonumber\\
 &\leq& \left\Vert S_h(s-r)(-A_h)^{\frac{1-\beta}{2}}\right\Vert^2_{L(H)}\left\Vert (-A_h)^{\frac{\beta-1}{2}} P_hQ^{\frac{1}{2}}\right\Vert^2_{\mathcal{L}_2(H)}\nonumber\\
 &\leq& C.
 \end{eqnarray}
 Using \assref{assumption5} and the smoothing properties of the semigroup, it holds that
 \begin{eqnarray}
 \label{mach6}
\mathbb{E}\left\Vert S_h(t_m-s)P_h\left(I^h_{k,s}-I^h_{k,k}\right)\right\Vert_{L(H)}
 &\leq&\left\Vert S_h(t_m-s)(-A_h)^{\frac{\eta}{2}}\right\Vert_{L(H)}\mathbb{E}\left\Vert(-A_h)^{-\frac{\eta}{2}}\left(I^h_{k,s}-I^h_{k,k}\right)\right\Vert_{L(H)}\nonumber\\
 &\leq& (t_m-s)^{-\frac{\eta}{2}}\mathbb{E}\left\Vert(-A_h)^{-\frac{\eta}{2}} \left(I^h_{k,s}-I^h_{k,k}\right)\right\Vert_{L(H)}.
 \end{eqnarray}
  From the definition of $I^h_{m,k,s}$ \eqref{Taylor2a}, by using   \lemref{distlemma}  we arrive at 
 \begin{eqnarray}
 \label{mach7}
 &&\Vert (-A_h)^{-\frac{\eta}{2}}\left(I^h_{k,s}-I^h_{k,k}\right)\Vert_{\mathcal{L}(H)}\nonumber\\
&\leq& \int_0^1\left\Vert (-A_h)^{-\frac{\eta}{2}}P_h\left( F'\left(X^h(t_k)+\lambda\left(X^h(s)-X^h(t_{k})\right)\right)-F'\left(X^h(t_{k})\right)\right)\right\Vert_{\mathcal{L}(H)}d\lambda\nonumber\\
&\leq& \int_0^1\left\Vert (-A_h)^{-\frac{\eta}{2}}\left( F'\left(X^h(t_{k})+\lambda\left(X^h(s)-X^h(t_{k})\right)\right)-F'\left(X^h(t_{k})\right)\right)\right\Vert_{\mathcal{L}(H)}d\lambda\nonumber\\
&\leq& C\int_0^1\lambda\Vert X^h(s)-X^h(t_{k})\Vert d\lambda=  C\Vert X^h(s)-X^h(t_{k})\Vert
 \end{eqnarray}
 Substituting \eqref{mach7} in \eqref{mach6} and  \thmref{timeregular} yields
 \begin{eqnarray}
 \label{mach8}
 \mathbb{E}\left\Vert S_h(t_m-s)P_h\left(I^h_{k,s}-I^h_{k,k}\right)\right\Vert^4_{L(H)}
 &\leq& C(t_m-s)^{-2\eta}\mathbb{E}\Vert X^h(s)-X^h(t_{k})\Vert^4\nonumber\\
 &\leq& C(t_m-s)^{-2\eta}(s-t_{k})^{\min(2, 2\beta)}.
 \end{eqnarray}
  Substituting \eqref{mach8} in \eqref{mach5} yields
 \begin{eqnarray}
 \label{mach9}
 \Vert III_{11}^{(32)}\Vert^2_{L^2(\Omega, H)}
 &\leq& C\sum_{k=0}^{m-1}\int_{t_{k}}^{t_{k+1}}\int_{t_{k}}^s(t_m-s)^{-\eta}(s-t_{k})^{\min(1, \beta)}drds\nonumber\\
 &\leq&  C\sum_{k=0}^{m-1}\int_{t_{k}}^{t_{k+1}}(t_m-s)^{-\eta}(s-t_{k})^{\min(2, 1+\beta)}ds\nonumber\\
 &\leq& C\Delta t^{\min(2, 1+\beta)}\sum_{k=0}^{m-1}\int_{t_{k}}^{t_{k+1}}(t_m-s)^{-\eta}\leq C\Delta t^{\min(2, 2\beta)}.
 \end{eqnarray}
 Substituting \eqref{mach9} and \eqref{Taylor7} in \eqref{Raphael1} yields
 \begin{eqnarray}
 \label{mach10}
 \Vert III_{11}^{(3)}\Vert_{L^2(\Omega, H)}\leq C\Delta t^{\frac{\beta}{2}}.
 \end{eqnarray}
 We split $III_{11}^{(4)}$ in two terms as follows
 \begin{eqnarray}
 \label{Raphael2}
 III_{11}^{(4)}&=&\sum_{k=0}^{m-1}\int_{t_k}^{t_{k+1}}S_h(t_m-s)I^h_{k,k}\int_{t_k}^s\int_{\chi}S_h(s-r)P_h\psi(z)\widetilde{N}(dz,dr)ds\nonumber\\
 &+&\sum_{k=0}^{m-1}\int_{t_k}^{t_{k+1}}S_h(t_m-s)\left(I^h_{k,s}-I^h_{k,k}\right)\int_{t_k}^s\int_{\chi}S_h(s-r)P_h\psi(z)\widetilde{N}(dz,dr)ds\nonumber\\
 &=:& III_{11}^{(41)}+III_{11}^{(42)}.
 \end{eqnarray}
Using the fact that the expectation of the cross-product terms vanishes, using H\"{o}lder inequality the It\^{o}-isometry, the smoothing properties of the semigroup, \lemref{distlemma}, it follows that
\begin{eqnarray}
\label{Taylor8}
\Vert III_{11}^{(41)}\Vert_{L^2(\Omega,H)}^2&=&\mathbb{E}
\left[\left\Vert\sum_{k=0}^{m-1}\int_{t_k}^{t_{k+1}}\left(\int_{t_k}^s\int_{\chi}S_h(t_m-s)I^h_{k,k}S_h(s-r)P_h\psi(z)\widetilde{N}(dz,dr)\right)ds\right\Vert^2\right]\nonumber\\
&=&\sum_{k=0}^{m-1}\mathbb{E}
\left\Vert\int_{t_k}^{t_{k+1}}\left(\int_{t_k}^s\int_{\chi}S_h(t_m-s)I^h_{k,k}S_h(s-r)P_h\psi(z)\widetilde{N}(dz,dr)\right)ds\right\Vert^2\nonumber\\
&\leq& \Delta t\sum_{k=0}^{m-1}\int_{t_k}^{t_{k+1}}\mathbb{E}\left\Vert\int_{t_k}^s\int_{\chi}S_h(s-t_k)I^h_{k,k}S_h(s-r)P_h\psi(z)\widetilde{N}(dz,dr)\right\Vert^2ds\nonumber\\
&\leq& \Delta t\sum_{k=0}^{m-1}\int_{t_k}^{t_{k+1}}\int_{t_k}^s\int_{\chi}\left\Vert S_h(t_m-s)I^h_{k,k}S_h(s-r)P_h\psi(z)\right\Vert^2\nu(dz)drds\nonumber\\
&\leq& \Delta t\sum_{k=0}^{m-1}\int_{t_k}^{t_{k+1}}\int_{t_k}^s\Vert S_h(t_m-s)I^h_{k,k}\Vert^2_{L(H)}\left\Vert S_h(s-r)(-A_h)^{\frac{1-\beta}{2}}\right\Vert^2_{L(H)}drds\nonumber\\
&&\times\int_{\chi}\left\Vert (-A_h)^{\frac{\beta-1}{2}}P_h\psi(z)\right\Vert^2\nu(dz)\nonumber\\
&\leq& C\Delta t\sum_{k=0}^{m-1}\int_{t_k}^{t_{k+1}}\int_{t_k}^s(s-r)^{\min(0,-1+\beta)}drds\leq C\Delta t^{\min(2,1+\beta)}.
\end{eqnarray}
Along the same lines as that of $III_{11}^{(32)}$ by using \eqref{Davis4}, one obtains 
\begin{eqnarray}
\label{Raphael3}
\Vert III_{11}^{(42)}\Vert^2_{L^2(\Omega, H)}\leq C\Delta t^{\min(2, 2\beta)}.
\end{eqnarray}
Substituting \eqref{Taylor8} and \eqref{Raphael3} in \eqref{Raphael2} yields
\begin{eqnarray}
\label{Raphael4}
\Vert III_{11}\Vert_{L^2(\Omega, H)}\leq C\Delta t^{\frac{\beta}{2}}.
\end{eqnarray}
Substituting \eqref{Raphael4}, \eqref{mach10}, \eqref{Taylor6} and \eqref{Taylor5} in \eqref{Taylor4} yields
\begin{eqnarray}
\label{lun18}
\Vert III_{11}\Vert_{L^2(\Omega, H)}\leq C\Delta t^{\beta/2}.
\end{eqnarray}
Substituting \eqref{lun18} and \eqref{lun4} in \eqref{lun3} yields 
\begin{eqnarray}
\label{lun19}
\Vert III_1\Vert_{L^2(\Omega, H)}\leq C\Delta t^{\beta/2}+C\Delta t\sum_{k=0}^{m-1}\Vert X^h(t_k)-X^h_k\Vert_{L^2(\Omega, H)}.
\end{eqnarray}
Let us now turn our attention to the the estimate of terms involving the noises. Let us start with $III_3$. Applying the It\^{o} isometry, the smoothing property of the semigroup, Lemmas \ref{distlemma} and \ref{sharpestimates} yields
\begin{eqnarray}
\label{lun20}
\Vert III_3\Vert^2_{L^2(\Omega, H)}&=&\sum_{k=0}^{m-1}\mathbb{E}\int_{t_k}^{t_{k+1}}\int_{\chi}\Vert S_h(t_m-s)\left(\mathbf{I}-S_h(s-t_k)\right)P_h\psi(z)\Vert^2\nu(dz)ds\nonumber\\
&\leq&\sum_{k=0}^{m-1}\int_{t_k}^{t_{k+1}}\int_{\chi}\Vert S_h(t_m-s)(-A_h)^{1/2}\Vert^2_{L(H)}\Vert\left(\mathbf{I}-S_h(s-t_k)\right)(-A_h)^{-\beta/2}\Vert^2_{L(H)}\nonumber\\
&&\times\Vert (-A_h)^{\frac{\beta-1}{2}}P_h\psi(z)\Vert^2\nu(dz)ds\nonumber\\
&\leq& C\sum_{k=0}^{m-1}\int_{t_k}^{t_{k+1}}\Vert S_h(t_m-s)(-A_h)^{1/2}\Vert^2_{L(H)}(s-t_k)^{\beta}ds\int_{\chi} \Vert (-A_h)^{\frac{\beta-1}{2}}P_h\psi(z)\Vert^2\nu(dz)\nonumber\\
&\leq& C\Delta t^{\beta}\sum_{k=0}^{m-1}\int_{t_k}^{t_{k+1}}\Vert S_h(t_m-s)(-A_h)^{1/2}\Vert^2_{L(H)}ds\nonumber\\
&\leq& C\Delta t^{\beta}\int_0^{t_m}\Vert S_h(t_m-s)(-A_h)^{1/2}\Vert^2_{L(H)}ds\nonumber\\
&\leq& C\Delta t^{\beta}.
\end{eqnarray}
Along the same lines as that of $III_2$, we obtain the following estimate
\begin{eqnarray}
\label{lun21}
\Vert III_2\Vert^2_{L^2(\Omega, H)}\leq C\Delta t^{\beta}.
\end{eqnarray}
Substituting \eqref{lun21}, \eqref{lun20} and \eqref{lun19} in \eqref{lun1} yields
\begin{eqnarray}
\label{lun22}
\Vert X^h(t_m)-X^h_m\Vert_{L^2(\Omega, H)}\leq C\Delta t^{\beta/2}+C\Delta t\sum_{k=0}^{m-1}\Vert X^h(t_k)-X^h_k\Vert_{L^2(\Omega, H)}.
\end{eqnarray}
Applying the discrete Gronwall lemma to \eqref{lun22} yields
\begin{eqnarray}
\label{lun23}
\Vert X^h(t_m)-X^h_m\Vert_{L^2(\Omega, H)}\leq C\Delta t^{\beta/2}.
\end{eqnarray}
Combining \eqref{lun23} and \thmref{spaceapproximation} completes the proof of \thmref{mainresult2}.

\subsection{Proofs  of \thmref{mainresult1}  and \thmref{mainresult2}  for semi implicit scheme}
\label{fulldiscretization2}
Let us recall that  the semi implicit Euler method was applied on SPDE \eqref{model} in \cite{Yang}  where the linear operator $A$ was assumed to be self adjoint, but the optimal convergence order was  
only sup-optimal of the form $\mathcal{O}(h^{1-\epsilon}+\Delta t^{1/2-\epsilon})$, for any $\epsilon>0$, small enough. 
Here,  we  upgrade  the proof   and we achieve optimal convergence rate $\mathcal{O}(h^2+\Delta t^{1/2})$ 
for SPDE \eqref{model} with multiplicative noises, and convergence order  $\mathcal{O}\left(h^2+\Delta t^{1-\epsilon}\right)$ for SPDE 
with additive noises  and not necessary self-adjoint operator.
Before sketching the proof of  Theorems \ref{mainresult1} and \ref{mainresult2} for our semi  implicit scheme,
 we need to recall  the following two lemmas some main ingredients from \cite{Antjd2}, 
 useful when dealing with not necessary self-adjoint operator. Note that the proof of \lemref{lemnonadditif1} can be found in \cite[Theorems 7.8 \& 7.8]{Vidar}, 
 but only for self-adjoint operator.  Note also that the proof of \lemref{lemnonadditif2} can be found in \cite[Lemma 4.4]{Raphael}, but only for self adjoint operator also. 

\begin{lemma}\cite{Antjd2}
\label{lemnonadditif1}
\begin{itemize}
\item[(i)] If $u\in\mathcal{D}((-A)^{\alpha/2})$, $0\leq \alpha\leq 2$, then it holds that
\begin{eqnarray}
\Vert\left(S_h(t_m)-S^m_{h,\Delta t}\right)P_hu\Vert\leq C\Delta t^{\alpha/2}\Vert u\Vert_{\alpha},\quad m=1,\cdots, M.
\end{eqnarray}
\item[(ii)] For $0\leq q\leq 1$ and $u\in H$, it holds that
\begin{eqnarray}
\Vert \left(S_h(t_m)-S^m_{h,\Delta t}\right)P_hu\Vert\leq C\Delta t^qt_m^{-q},\quad m=1,\cdots, M.
\end{eqnarray}
\end{itemize}
\end{lemma}
\begin{lemma}\cite{Antjd2}
\label{lemnonadditif2}
Let $0\leq \rho\leq 1$ and $\epsilon>0$ an arbitrarily positive number, small enough.
\begin{itemize}
\item[(i)] For all $u\in\mathcal{D}((-A)^{-\rho})$, we have 
\begin{eqnarray}
\left\Vert\sum_{j=1}^m\int_{t_{j-1}}^{t_j}\left(S^j_{h,\Delta t}-S_h(s)\right)P_huds\right\Vert\leq C\Delta t^{\frac{2-\rho}{2}-\epsilon}\Vert u\Vert_{-\rho}.
\end{eqnarray}
\item[(ii)] For any $\gamma\in[0, 2]$ and $u\in\mathcal{D}((-A)^{\frac{\gamma-1}{2}}$, the following estimate holds
\begin{eqnarray}
\left(\sum_{j=1}^m\int_{t_{j-1}}^{t_j}\left\Vert\left(S^j_{h, \Delta t}-S_h(s)\right)P_hu\right\Vert^2ds\right)^{1/2}\leq C\Delta t^{\gamma/2-\epsilon}\Vert u\Vert_{\gamma-1}.
\end{eqnarray}
\end{itemize}
\end{lemma}
\begin{proof}
\textit{[of Theorems \ref{mainresult1} and \ref{mainresult2} for semi  implicit scheme]}

The proof of \thmref{mainresult1} for $\beta\in[0,1)$ follows the sames lines as in \cite[Theorem 6.1]{Yang} by using \lemref{lemnonadditif1}. Note that the proof in the case $\beta\in[1,2]$ makes use of \lemref{sharpestimates}. 
The proof of \thmref{mainresult2} follows the lines of the proof of the main result in \cite{Xiaojie2}, by using \lemref{lemnonadditif2} and some results in \cite{Antjd2}, appropriate when dealing with not necessarily self-adjoint operator. 
Note that the term involving the Poisson measure is handled  in the same way as in the proofs of Theorems  \ref{mainresult1} and \ref{mainresult2} for exponential scheme.
\end{proof}
\section{Numerical experiments}
\label{numericalexperiments}

\subsection{Additive Gaussian and Poison measure noises}
Here, we consider  the stochastic advection-diffusion-reaction SPDE with additive  noises in two dimensions on the domain $\Lambda=[0,L_l]\times[0,L_l]$. 
\begin{eqnarray}
\label{reactiondif2}
dX&=&\left[\nabla \cdot (\mathbf{D}\nabla X)-\nabla \cdot(\mathbf{q}X) +\dfrac{X}{1+X^2}\right]dt+dW+z_0d\widetilde{N}, 
\end{eqnarray}
with $D=10^{-1} \mathbf{I}_2$, where $\mathbf{I}_2$ is $2 \times 2$ matrice.  The Poisson counting process has a constant  rate $\lambda=0.2$ with finite measure.
In the abstract setting \eqref{model1}, our linear operator $A$ is given by 
\begin{eqnarray}
 A=\nabla \cdot  \mathbf{D}\nabla (.) -  \nabla \cdot \mathbf{q}(.).
\end{eqnarray}
The Darcy velocity  $\mathbf{q}$ is obtained as in \cite{Antonio1}  and  to deal with high  P\'{e}clet flows.
Let us write $A=A_s+A_{ns}$, where $A_s$ is the self-adjoint part of $A$ and $A_{ns}$ its  non self-adjoint part. An easy computation shows that the eigenfunctions $\{e_{i,j}\}_{i, j\geq 0}=\{e_i^{(1)}\otimes e_j^{(2)}\}_{i,j\geq 0}$ with the corresponding eigenvalues $\{\lambda_{i,j}\}=\{(\lambda_i^{(1)})^2+(\lambda_j^{(2)})^2\}$ of $-A_s$ are given by
\begin{eqnarray}
\label{eigen1}
e_i^{(l)}(x)=\frac{1}{\sqrt{L_l}},\quad \lambda_i^{(l)}=\frac{iD\pi}{L_l},\quad l=1,2,\quad x\in \Lambda,\quad i\in\mathbb{N}.
\end{eqnarray}
Note $\mathcal{D}(-A)=\mathcal{D}(-A_s)$ with equivalent norms (see e.g. \cite{Suzuki}), so by \cite[(3.3)]{Lions} we have $\mathcal{D}((-A)^{\alpha})=\mathcal{D}((-A_s)^{\alpha})$, $0\leq \alpha\leq 1$ with equivalence of norms. Therefore 
\begin{eqnarray}
\label{eigen3}
\Vert (-A)^{\frac{\beta-1}{2}}Q^{\frac{1}{2}}\Vert_{\mathcal{L}_2(H)}\leq C\Vert (-A_s)^{\frac{\beta-1}{2}}Q^{\frac{1}{2}}\Vert_{\mathcal{L}_2(H)}.
\end{eqnarray}
In our example \eqref{reactiondif2}, we consider the covariance operator $Q$ to have the same eigenfunctions as $-A_s$ with the eigenvalues given by
\begin{eqnarray}
\label{eigen2}
q_{i,j}=\frac{1}{(i^2+j^2)^{\beta+\delta}},\quad 0\leq \beta\leq 2, \quad \delta>0,\quad \text{small enough}.
\end{eqnarray}
Therefore, using \eqref{eigen3} one can easily check that
\begin{eqnarray}
\Vert (-A)^{\frac{\beta-1}{2}}Q^{\frac{1}{2}}\Vert^2_{\mathcal{L}_2(H)}<\infty.
\end{eqnarray}
Using the equivalence of norms $\Vert (-A)^{\alpha}.\Vert$ and $\Vert (-A_s)^{\alpha}.\Vert$, $0\leq \alpha\leq 1$, writing $z_0\in H$ of \eqref{reactiondif2} in eigenbasis $(e_{i,j})_{i, j\geq 0}$ of the operator $-A_s$ and using Parseval's identity yields
\begin{eqnarray}
\label{eigen4}
\Vert (-A)^{\frac{\beta-1}{2}}z_0\Vert^2&\leq& C\Vert (-A_s)^{\frac{\beta-1}{2}}z_0\Vert^2=C\left\Vert\sum_{(i,j)\in\mathbb{N}^2}(-A_s)^{\frac{\beta-1}{2}}\left\langle z_0, e_{i,j}\right\rangle_H e_{i,j}\right\Vert^2\nonumber\\
&=&C\sum_{(i,j)\in\mathbb{N}}\lambda_{i,j}^{\beta-1}\langle z_0, e_{i,j}\rangle_H^2=\frac{CD^2\pi^2}{L_l^2}\sum_{(i,j)\in\mathbb{N}^2}(i^2+j^2)^{\beta-1}\langle z_0, e_{i,j}\rangle_H^2.
\end{eqnarray}
Taking $z_0$ such that 
\begin{eqnarray}
\label{eigen5}
\langle z_0, e_{i,j}\rangle=\frac{1}{(i^2+j^2)^{\frac{\beta+\gamma}{2}}},\quad 0\leq\gamma\leq 2,\quad \gamma>0\quad \text{small enough},
\end{eqnarray}
it follows from \eqref{eigen4} that
\begin{eqnarray}
\label{eigen6}
\Vert (-A)^{\frac{\beta-1}{2}}z_0\Vert^2\leq \frac{CD^2\pi^2}{L_l^2}\sum_{(i,j)\in\mathbb{N}^2}(i^2+j^2)^{-(1+\gamma)}<\infty,\quad \gamma>0.
\end{eqnarray}
Therefore for
\begin{eqnarray}
z_0=\sum_{(i,j)\in\mathbb{N}^2}(i^2+j^2)^{-\frac{(\beta+\gamma)}{2}}e_{i,j},\quad \gamma>0
\end{eqnarray}
 \assref{assumption5} is fulfilled. 
 The drift function is given by $F(u)=\frac{u}{1+u^2}$, $u\in H$; and one can easily check that \assref{assumption2}  is  satisfied. Here, the Nemytskii operator defined in \eqref{nemystskii} is given by $f(x,y)=\frac{y}{1+y^2}$ for all $(x, y)\in\Lambda\times\mathbb{R}$.  One can easily check that
 \begin{eqnarray}
 \label{borne1}
 \vert f(x,y)\vert\leq C(1+\vert y\vert),\quad \frac{\partial f}{\partial y}(x,y)=\frac{1-y^2}{(1+y^2)^2},\quad \left\vert\frac{\partial f}{\partial y}(x,y)\right\vert\leq C,\quad (x,y)\in\Lambda\times\mathbb{R}.
 \end{eqnarray}
 From \eqref{nemystskii}, it follows that
 \begin{eqnarray}
 \label{borne3}
 F'(v)(u)(x)&=&\frac{\partial f}{\partial y}(x, v(x)).u(x),\quad x\in\Lambda,\quad u,v\in H.
 \end{eqnarray}
 It follows from \eqref{borne3} by using \eqref{borne1} that
\begin{eqnarray}
\label{borne5}
\Vert F'(v)u\Vert\leq C\Vert u\Vert,\quad u, v\in H.
\end{eqnarray}
One can also easily check that 
\begin{eqnarray*}
\Vert (-A)^{-\frac{\eta}{2}}\left(F'(u)-F'(v)\right)\Vert_{L(H)}\leq C\Vert u-v\Vert,\quad u,v\in H.
\end{eqnarray*}
 Therefore \assref{assumption6} is fulfilled.

   We use in our simulations $\beta=2$ and $\gamma=10^{-3}$ in \eqref{eigen5}.
 We discretize in space using finite volume method (viewed as the finite element method as in \cite{Antonio1}) in rectangular grid of sizes $\Delta x=\Delta y= 1/300$. 
\begin{figure}[!ht]
 \begin{center}
    \label{FIGI}
    \includegraphics[width=0.47\textwidth]{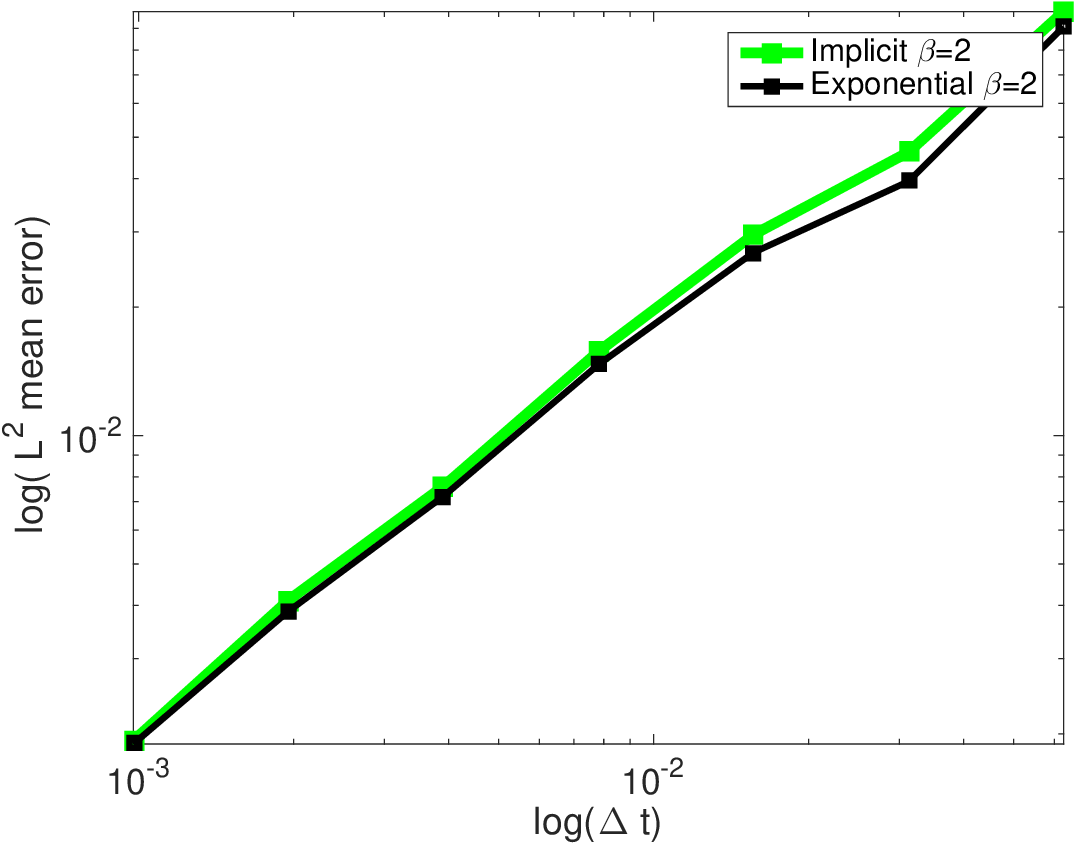}
  \caption{(a) Convergence in the root mean square $L^{2}$ norm at $T=1$ as a
    function of $\Dt$. We show convergence for noise  where 
    $\beta=2$, and  $\delta=10^{-3}=\gamma$ in relation \eqref{eigen2}. We have  used here 30 realizations.  The order of convergence is $0.95$ for exponential scheme and $0.92$ for semi implicit scheme. 
    Note that the intial solution is $X_0=0$
    } 
  \label{FIGI} 
  \end{center}
  \end{figure}
   We take $L_l=1$. The convergence graphs for the exponential and semi implicit scheme are presented in  \figref{FIGI}. The order of convergence in time is $0.95$ and $0.92$ respectively.  
  These orders are close to the theoretical  orders in \thmref{mainresult2}.

\subsection{Multiplicative Gaussian and Poisson measure noises}
Let us consider  the stochastic advection-diffusion-reaction SPDE with multiplicative noises in two dimensions on the domain $\Lambda=[0,2]\times[0,3]$. 
\begin{eqnarray}
\label{reactiondif1}
dX&=&\left[\nabla \cdot (\mathbf{D}\nabla X)-\nabla \cdot(\mathbf{q}X) +\dfrac{X}{1+X^2}\right]dt+XdW+Xd\widetilde{N}.\\
\mathbf{D}&=&\left( \begin{array}{cc}
             10^{-1}&0\\
             0& 10^{-2}
             \end{array}\right)
\end{eqnarray}
with mixed Neumann-Dirichlet boundary conditions. The Dirichlet boundary condition is $X=1$ at $x=0$. 
 We use the homogeneous Neumann boundary conditions elsewhere. The Darcy velocity  $\mathbf{q}$ is obtained as in \cite{Antonio1}  and  to deal with high  P\'{e}clet flows.
We discretize in space using finite volume method (viewed as the finite element method as in \cite{Antonio1}) in rectangular grid of sizes $\Delta x=\Delta y= 1/300$. 
The noises are the same as for additive noise described in the previous section.  The Poisson counting process has a constant  rate $\lambda=0.2$.
\begin{figure}[!ht]
    \label{FIGII}
    \includegraphics[width=0.47\textwidth]{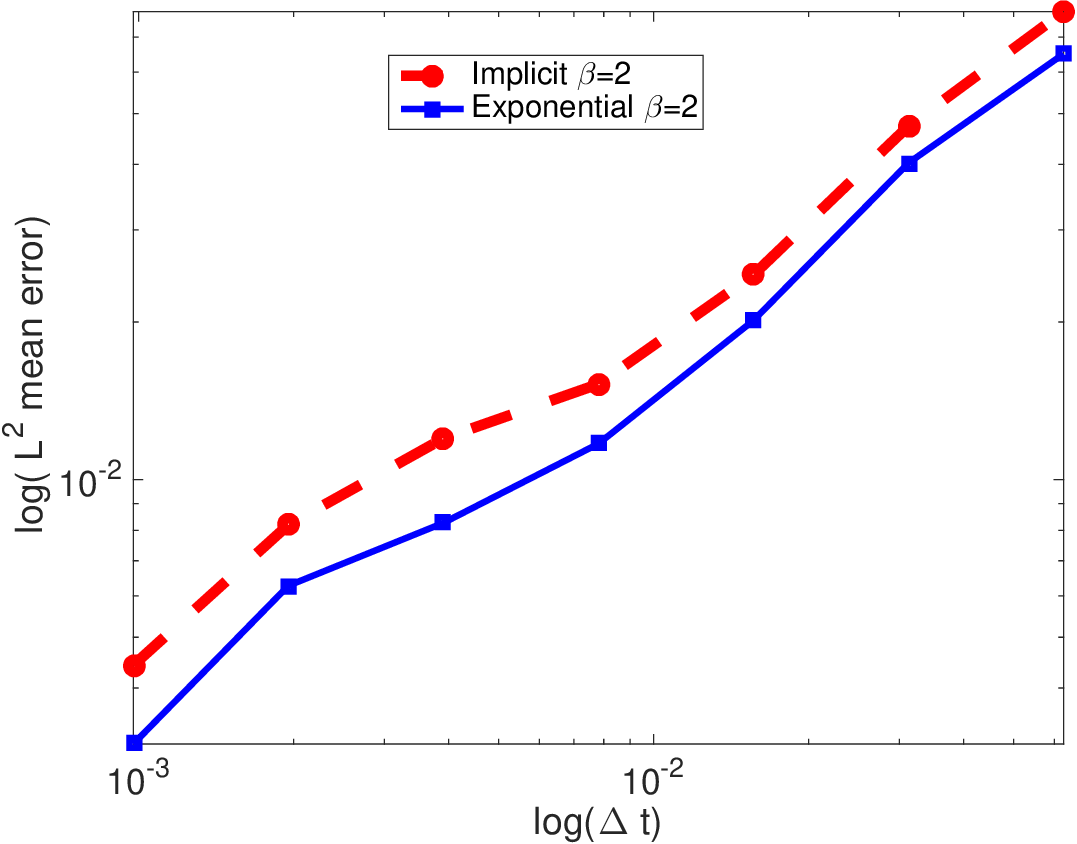}
  \caption{(a) Convergence in the root mean square $L^{2}$ norm at $T=1$ as a
    function of $\Dt$. We show convergence for noise  where 
    $\beta=2$, and  $\delta=10^{-3}$ in relation \eqref{eigen2}. We have  used here 30 realizations.  The order of convergence is $0.58$ for exponential scheme and $0.55$ for semi implicit scheme.
    Note that the initial solution is $X_0=0$ 
    } 
  \label{FIGII} 
  \end{figure}
 The convergence graphs for the exponential and semi implicit scheme are presented in \figref{FIGII}. The order of convergence in time is $0.58$ and $0.55$ respectively.  
  These orders are close to the theoretical  orders in \thmref{mainresult1}.

\section*{Acknowledgements}
The authors thank the reviewers for their careful reading and comments that help to improve the manuscript.  Part of this work was done when J. D. Mukam visited the Western Norway University of Applied Sciences. The visit  was supported by the "InPro TU Chemnitz program". J. D. Mukam thanks the "InPro TUC program" for making this stay possible.

\end{document}